\documentclass[preprint,sort&compress,final,10pt]{elsarticle}

\usepackage{amsmath}
\usepackage{amssymb}
\usepackage{graphicx}
\usepackage{latexsym}
\usepackage{epstopdf}
\usepackage{natbib}
\usepackage{color}
\usepackage{hyperref}

\newtheorem{theorem}{Theorem}[section]

\newtheorem{lemma}[theorem]{Lemma}
\newtheorem{proposition}[theorem]{Proposition}

\numberwithin{equation}{section}

\def\Proof{\noindent{\bf Proof.}~}
\def\qed{\hfill$\square$\smallskip}
\def\dint{\displaystyle\int}

\def\dsum{\displaystyle\sum}
\def\dint{\displaystyle\int}

\def\Im{\mathrm{Im}}

\journal{\empty}
\date{}

\begin{document}

\begin{frontmatter}

\title{Boundedness of solutions in impulsive Duffing equations with polynomial potentials and $C^{1}$ time dependent coefficients}

\author[au1]{Yanmin Niu}

\ead[au1]{nym1216@163.com}

\author[au1]{Xiong Li\footnote{ Partially supported by the NSFC (11571041)and the Fundamental Research Funds for the Central Universities. Corresponding author.}}

\address[au1]{School of Mathematical Sciences, Beijing Normal University, Beijing 100875, P.R. China.}

\ead[au1]{xli@bnu.edu.cn}

\begin{abstract}
In this paper, we are concerned with the impulsive Duffing equation
$$
x''+x^{2n+1}+\sum_{i=0}^{2n}x^{i}p_{i}(t)=0,\ t\neq t_{j},
$$
with impulsive effects $x(t_{j}+)=x(t_{j}-),\ x'(t_{j}+)=-x'(t_{j}-),\ j=\pm1,\pm2,\cdots$, where the time dependent coefficients $p_i(t)\in C^1(\mathbb{S}^1)\ (n+1\leq i\leq 2n)$ and  $p_i(t)\in C^0(\mathbb{S}^1)\ (0\leq i\leq n)$ with $\mathbb{S}^1=\mathbb{R}/\mathbb{Z}$. If impulsive times are 1-periodic and $t_{2}-t_{1}\neq\frac{1}{2}$ for $0< t_{1}<t_{2}<1$, basing on a so-called large twist theorem recently established by X. Li, B. Liu and Y. Sun in \cite{XLi}, we find large invariant curves diffeomorphism to circles surrounding the origin and going to infinity, which confines the solutions in its interior and therefore leads to the boundedness of these solutions. Meanwhile, it turns out that the solutions starting at $t=0$ on the invariant curves are quasiperiodic.
\end{abstract}


\end{frontmatter}

\section{Introduction and main results}
In this paper, we discuss the impulsive Duffing equations
\begin{equation}\label{maineq}
\left\{\begin{array}{ll}
 x''+x^{2n+1}+\displaystyle\sum_{i=0}^{2n}x^{i}p_{i}(t)=0,\quad t\neq t_{j};\\[0.2cm]
 \Delta x|_{t=t_{j}}=0,\\[0.2cm]
 \Delta x'|_{t=t_{j}}=-2x'(t_{j}-),\quad\ \ \ \ \ j=\pm1,\pm2,\cdots,\\
 \end{array}
 \right.
 \end{equation}
where $0< t_{1}<t_{2}<1$ satisfying $t_{2}-t_{1}\neq\frac{1}{2}$, $n\geq1$, $\Delta x|_{t=t_{j}}=x\left(t_{j}+\right)-x\left(t_{j}-\right)$ and $\Delta x'|_{t=t_{j}}=x'\left(t_{j}+\right)-x'\left(t_{j}-\right)$. In addition, assume that the impulsive time is $1$-period, that is, $t_{j+2}=t_{j}+1$ for $j=\pm1,\pm2,\cdots$, and $p_i(t)\in C^1(\mathbb{S}^1) \,(n+1\leq i\leq 2n)$ and  $p_i(t)\in C^0(\mathbb{S}^1)\, (0\leq i\leq n)$ with $\mathbb{S}^1=\mathbb{R}/\mathbb{Z}$.
The behaviour of a time dependent nonlinear differential equation
\begin{equation}
x''+f(t,x)=0,
\end{equation}
$f$ being periodic in $t$, can be very intricate. For example, there are equations having unbounded solutions but with infinitely many zeroes and with unbounded solutions having randomly prescribed numbers of zeroes and also periodic solutions, see \cite{alekseev}, \cite{moser} and \cite{sitnikov}. Besides such unbounded phenomena, many authors take notice on looking for conditions on the nonlinearity, which allow to conclude that all the solutions of the equation are bounded, see \cite{Dieckerhoff}, \cite{Laederich}, \cite{Yuan95}-\cite{Yuan00}.
It is well known that every solution of the Duffing equation
\begin{equation}
x''+x^{3}=p(t),
\end{equation}
$p(t+1)=p(t)$ being continuous, is bounded. This result, prompted by questions of J. E. Littlewood in \cite{littlewood}, is due to G. R. Morris \cite{morris}.
In \cite{Dieckerhoff}, R. Dieckerhoff and E. Zehnder further study the Duffing equations
\begin{equation}\label{duffing}
x''+x^{2n+1}+\displaystyle\sum_{j=0}^{2n}x^{j}p_{j}(t)=0,\ \ \ \ n\geq1,
\end{equation}
where $p_j(t)$ are
required to be sufficiently smooth to construct a series change of variables to
transform this equation into a nearly integrable system for large energies. By twist theorem, they find large invariant curves and therefore obtain the boundedness of solutions and quasiperiodic solutions.

Among these studies for time dependent Duffing equations, an interesting problem (proposed by Diecherhoff and Zehnder) is whether or not the boundedness of all solutions depends on the smoothness of the coefficients. In fact in \cite{Dieckerhoff}, the author remark the smoothness on $p_i(t)$ depends on the index $i$. In \cite{Laederich}, Laederich and Levi weakened the smoothness requirement on $p_i(t)$ to $C^{5+\epsilon}$. By modifying the proofs in \cite{Dieckerhoff} and using some approximation techniques, B. Liu \cite{Liu} obtained the same result for
$$
x'' + x^{2n+1}+p_1(t)x+p_0(t) = 0,
$$
where the periodic functions $p_1,\ p_0$ are only required to be continuous, which shows that the boundedness of all solutions does not depend on the smoothness of coefficients of lower order terms. Later, Yuan \cite{Yuan95}, \cite{Yuan98}, \cite{Yuan00}
proved that all solutions of (\ref{duffing}) are bounded if $p_i(t)\in C^2\,(n+1 \leq i \leq 2n)$ and $p_i(t) \in C^1 \,(0 \leq i \leq n)$. Recently in \cite{XLi}, X. Li, B. Liu and Y. Sun show the boundedness of (\ref{duffing}) under more weaker conditions $p_i(t)\in C^1(\mathbb{S}^1)\ (n+1\leq i\leq 2n)$ and  $p_i(t)\in C^0(\mathbb{S}^1)\ (0\leq i\leq n)$ with $\mathbb{S}^1=\mathbb{R}/\mathbb{Z}$. Of course, there are papers explain the vital importance of smoothness of higher order terms' coefficients for the boundedness of Duffing equations' solutions from the reverse side. In \cite{Levi}, Levi and You proved that the equation
$$
x'' + x^{2n+1} + p(t)x^{2i+1} = 0
$$
with a special discontinuous coefficient $p(t) = K^{[t]mod2}, 0 < K < 1, 2n + 1 >
2i + 1 \geq  n + 2$, possesses an oscillatory unbounded solution. In 2000, Wang \cite{Wang00}
constructed a continuous periodic function $p(t)$ such that the corresponding
equation
$$
x'' + x^{2n+1} + p(t)x^{i} = 0
$$
possesses a solution which escapes to infinity in finite time, where $n \geq 2$ and $2n + 1 > i\geq n + 2$. Both the two examples show that the boundedness of solutions is linked with the smoothness of coefficients of higher order terms.

Our aim is to extend the result in \cite{Dieckerhoff} under weaker smoothness conditions of coefficients to the case with impulsive effects. Recently, as impulsive equations widely arise in applied mathematics, they attract a lot of attentions and many authors study the general properties of them in \cite{Bainov93}, \cite{Lak}, along with the existence of periodic solutions via fixed point theory in \cite{Nieto97}, \cite{Nieto02}, topological degree theory in \cite{Dong01}, \cite{Qian05}, and variational method in \cite{Nieto09}, \cite{Sun11}. However, different from the extensive study for second order differential equations without impulsive terms, there are only a few results on boundedness of solutions and the existence of quasiperiodic solutions for  second order impulsive differential equations.

The existence of impulses, even the simplest impulsive functions, may cause complicated dynamic phenomenons and bring difficulties to study. The behavior of solutions with impulsive effects may have great differences compared with solutions without impulsive effects. Choosing different impulsive functions may has different effects on the solutions. The construction of impulsive functions in \eqref{maineq} is inspired by linear oscillator equation with obstacle, which has the form
\begin{equation}\label{oscillator}
\left\{\begin{array}{ll}
 x''+ax=f(t),\\[0.2cm]
 x(t)\geq 0\\[0.2cm]
 x(t_{0})=0\ \Rightarrow x'(t_{0}+)=-x'(t_{0}-),\\
 \end{array}
 \right.
 \end{equation}
where $a>0$ and $f\in C(\mathbb{S}^1)$. This equation has actual physical backgrounds in classical billiard in electric and magnetic fields, three body problem and so on, see \cite{pboyland}, \cite{corbera}, \cite{hlamba}. Equation \eqref{oscillator} can be thought of as the model of the motion of a particle which is attached to a spring that pushes the particle against a rigid wall situated at $x=0$. At this barrier the particle bounces elastically. At the time $t_{0}$ satisfying $x(t_{0})=0$, the change $x'(t_{0}+)=-x'(t_{0}-)$ is similar with the impulsive function in \eqref{maineq}. In both cases, the motions of solutions are actually continuous and just the rates of velocity change signs. However, in \eqref{maineq} impulsive times are pre-given while in \eqref{oscillator} the positions where solutions meet the obstacle and change their orientations are fixed. From this point, the two cases have essential difference.

Another reason we choose these impulsive functions is that it reduces the difficulty causing by impulses. It is obvious that the solutions of autonomous Duffing equation are a family of closed curves in the phase planer. Due to the existence of impulses, the closed curves cannot be preserved and even be broken in general. For this point, the impulsive functions in \eqref{maineq} guarantee the motions of trajectories still on the closed curves, which makes the study about such equations accessible.

In \cite{niu}, Y. Niu and X. Li discuss the periodic solutions of semilinear impulsive differential equation
 \begin{equation}\label{niu}
\left\{\begin{array}{ll}
 x''+g(x)=p(t),\quad t\neq t_{j};\\[0.2cm]
 \Delta x|_{t=t_{j}}=0,\\[0.2cm]
 \Delta x'|_{t=t_{j}}=-2x'\left(t_{j}-\right),\quad\ \ \ \ \ j=\pm1,\pm2,\cdots,\\
 \end{array}
 \right.
 \end{equation}
where $0<t_{1}<2\pi$, $\Delta x|_{t=t_{j}}=x\left(t_{j}+\right)-x\left(t_{j}-\right)$ and
$\Delta x'|_{t=t_{j}}=x'\left(t_{j}+\right)-x'\left(t_{j}-\right)$, $g(x),\ p(t)\in C(\mathbb{R},\mathbb{R})$ and $p(t)$ is $2\pi$-periodic. In addition, they assume that the impulsive time is $2\pi$-periodic, that is,
$t_{j+1}=t_{j}+2\pi$ for $j=\pm1,\pm2,\cdots$. By the Poincar\'{e}-Birkhoff twist theorem, they obtain finitely many $2\pi$-periodic solutions.

In this paper, we choose the same impulsive conditions as in \eqref{niu}, but the nonlinearity is superlinear and we use a so-called large twist theorem through a serious of symplectic transformations to obtain the boundedness of solutions and quasiperiodic solutions. The main results are as follows.

\begin{theorem}\label{boundedness}
If $p_i(t)\in C^1(\mathbb{S}^1) \,(n+1\leq i\leq 2n)$ and  $p_i(t)\in C^0(\mathbb{S}^1)\, (0\leq i\leq n)$ and $t_j\ (j=\pm1,\pm2,\cdots)$ are defined same as in \eqref{maineq} with $t_{2}-t_{1}\neq\frac{1}{2}$, then every solution $x(t)$ of \eqref{maineq} is bounded, i.e. it exists for all $t\in\mathbb{R}$ and
$$
\sup_{\mathbb{R}/\{t_{j}\}}(|x(t)|+|x'(t)|)+\sup_{j=\pm1,\pm2,\cdots}(|x(t_{j})|+|x'(t_{j})|)<\infty.
$$
\end{theorem}

\begin{theorem}\label{quasiperiodic solution}
Under the same assumptions as in Theorem \ref{boundedness}, there is a large $\omega^{\ast}>0$ such that for every irrational number $\omega>\omega^{\ast}$ satisfying
$$
\big|\omega-\frac{p}{q}\big|\geq c|q|^{-2-\beta}
$$
for all integers $p$ and $q\neq 0$ with two constants $\beta>0$ and $c>0$, there is a quasiperiodic solution of \eqref{maineq} having frequencies $(\omega,1)$; i.e. there is a function $F(\theta_{1},\theta_{2})$ periodic of period $1$, such that
$$
x(t)=F(\omega t, t)
$$
are solutions of the equation.
\end{theorem}

\begin{theorem}\label{hamonic period solutions}
Under the same assumptions as in Theorem \ref{boundedness}, for every integer $m\geq1$, there are infinitely many periodic solutions of \eqref{maineq} having minimal period $m$.
\end{theorem}

\section{Action and angle-variables}
Dropping the time dependent terms and impulse effects, we first consider equation
$$
 x''+x^{2n+1}=0,
$$
which can be written as an equivalent system
\begin{equation}\label{autonomou two eq}
\left\{\begin{array}{ll}
 x'=y,\ \ \\[0.2cm]
 y'=-x^{2n+1}.
\end{array}
\right.
\end{equation}
This is a time independent Hamiltonian system on $\mathbb{R}^{2}$:
\begin{equation}\label{H system}
\begin{array}{ll}
X_{h_{0}}:\left\{\quad\begin{array}{ll}
x'=\frac{\partial}{\partial y}h_{0}(x,y), \\[0.2cm]
y'=-\frac{\partial}{\partial x}h_{0}(x,y),\\[0.2cm]
\end{array}
\right.
\end{array}
 \end{equation}
where
$$
h_{0}(x,y)=\frac{1}{2}y^{2}+\frac{1}{2(n+1)}x^{2(n+1)}.
$$ Clearly $h_{0}>0$ on $\mathbb{R}^{2}$ except at the only equilibrium point $(x,y)=(0,0)$ where $h_{0}=0$. All the solutions of \eqref{H system} are periodic, the periods tending to zero as $h_{0}=E$ tends to infinity. We take $(x^{\ast}(t),y^{\ast}(t))$ for a solution of $X_{h_{0}}$ system having the initial conditions $(x^{\ast}(0),y^{\ast}(0))=(1,0)$. It has the energy $h_{0}(x^{\ast}(t),y^{\ast}(t))=E=1/2(n+1)$. Let $T^{\ast}>0$ be its minimal period and introduce the functions $C$ and $S$ by
$$
(C(t),S(t))=(x^{\ast}(t),y^{\ast}(t))
$$
These analytic functions satisfy
$$
\begin{array}{ll}
(i)\ C(t)=C(t+T^{\ast}),\ S(t)=S(t+T^{\ast})\ \mbox{and}\ C(0)=1,\ S(0)=0,\\[0.2cm]
(ii)\ C'(t)=S(t),\ S'(t)=-C(t)^{2n+1},\\[0.2cm]
(iii)\ (n+1)S(t)^{2}+C(t)^{2(n+1)}=1,\\[0.2cm]
(iv)\ C(-t)=C(t),\ \mbox{and}\ S(-t)=-S(t).
\end{array}
$$
The action and angle variables are now defined by the map $\varphi:\mathbb{R}\times \mathbb{S}^{1}\rightarrow\mathbb{R}^{2}\setminus0$, where $(x,y)=\varphi(\lambda,\theta)$ with $\lambda>0$ and with $\theta\ (mod\ 1)$ is given by the formulae:
\begin{equation}\label{symplectic tran}
\begin{array}{ll}
\varphi:\quad \begin{array}{ll}
 x=c^{\alpha}\lambda^{\alpha}C(\theta T^{\ast}),\ \ \ \ \\[0.2cm]
 y=c^{\beta}\lambda^{\beta}S(\theta T^{\ast}),\ \ \\[0.2cm]
\end{array}
 \end{array}
 \end{equation}
where $\alpha=\frac{1}{n+2}$, $\beta=1-\alpha$ and $c=1/(\alpha T^{\ast})$.

We claim that $\varphi$ is a symplectic diffeomorphism from $\mathbb{R}^{+}\times \mathbb{S}^{1}$ onto $\mathbb{R}^{2}\setminus0$. Indeed, for the Jacobian $\Delta$ of $\varphi$ one finds by $(ii)$ and $(iii)$ $|\Delta|=1$, so that $\varphi$ is measure preserving. Moreover since $(C,S)$ is a solution of a differential equation having $T^{\ast}$ as minimal period one concludes that $\varphi$ is one to one and onto, which proves the claim.

As for \eqref{maineq} without time dependent terms, it can be written in the form of a Hamiltonian system
\begin{equation}\label{Xh}
\begin{array}{ll}
 \left\{ \begin{array}{ll}
x'=\frac{\partial}{\partial y}h(x,y), \\[0.2cm]
y'=-\frac{\partial}{\partial x}h(x,y),\ \ \ t\neq t_{j}, \\[0.2cm]
x(t_{j}+)=x(t_{j}-),\ \ \\[0.2cm]
y(t_{j}+)=-y(t_{j}-),\ \ \ j=\pm1,\pm2,\cdots,
 \end{array}
 \right.
 \end{array}
 \end{equation}
where
$$
h(x,y)=\frac{1}{2}y^{2}+\frac{1}{2(n+1)}x^{2(n+1)}.
$$

The trajectories of \eqref{Xh} in $(x,y)$ planer are closed curves which are symmetric with respect to the $x$-axis and $y$-axis. Under the influence of the special impulses in this equation, motions of the solutions cannot escape the closed curve. Moreover the solutions are smooth while the derivatives change signs at the impulsive times. In the new coordinate $(\lambda,\theta)$, action $\lambda$ is the area encircled by the solution orbit, whereas $\theta$ is the angle with the $x$-axis. Then at $t=t_{j}$,
$$
\lambda(t_{j}+)=\lambda(t_{j}-),\ \ \ \ \theta(t_{j}+)=-\theta(t_{j}-),\ \ j=\pm1,\pm2,\cdots.
$$
Under the symplectic transformation $\psi$, the Hamiltonian function of \eqref{Xh} becomes
$$
h\circ\psi(\lambda,\theta)=d\lambda^{2\beta}\triangleq h^{\ast}(\lambda),
$$
while $X_{h}$ system becomes very simple:
\begin{equation}\label{H impulse system under action angle}
\begin{array}{ll}
X_{h^{\ast}}:\left\{
\begin{array}{ll}
\theta'=\frac{\partial}{\partial \lambda}h^{\ast}=2\beta d\lambda^{2\beta-1},\ \ \\[0.2cm]
\lambda'=-\frac{\partial}{\partial \theta}h^{\ast}=0,\quad \ t\neq t_{j};\\[0.2cm]
\lambda(t_{j}+)=\lambda(t_{j}-),\\[0.2cm]
\theta(t_{j}+)=-\theta(t_{j}-),\ \ \ j=\pm1,\pm2,\cdots.
\end{array}
\right.
\end{array}
\end{equation}

The full equation \eqref{maineq} has the form of
\begin{equation}\label{Xhfull}
\begin{array}{ll}
X_{h}: \left\{
\begin{array}{ll}
 x'=y,\ \ \\[0.2cm]
 y'=-x^{2n+1}-\sum_{j=0}^{2n}x^{j}p_{j}(t),\quad t\neq t_{j};\\[0.2cm]
x(t_{j}+)=x(t_{j}-),\\[0.2cm]
y(t_{j}+)=-y(t_{j}-),\ \ \ j=\pm1,\pm2,\cdots,
 \end{array}
 \right.
 \end{array}
 \end{equation}
with the Hamiltonian function:
$$
h(x,y,t)=\frac{1}{2}y^{2}+\frac{1}{2(n+1)}x^{2(n+1)}+\sum_{i=0}^{2n}\frac{p_i(t)}{i+1}x^{i+1},\quad t\neq t_{j},\ \ j=\pm1,\pm2,\cdots.
$$
Under the symplectic transformation $\psi$, $h$ becomes:
\begin{equation}\label{H}
H(\lambda,\theta,t)\triangleq h(\psi(\lambda,\theta),t)=d\lambda^{2\beta}+h_1(\lambda,\theta,t) + h_2(\lambda,\theta,t),\quad t\neq t_{j},\ \ j=\pm1,\pm2,\cdots,
\end{equation}
where
$$
h_1(\lambda,\theta,t)=\sum_{i=n+1}^{2n}a_i\lambda^{(i+1)\alpha}C(\theta T^*)^{i+1}p_i(t),
$$
$$
h_2(\lambda,\theta,t)=\sum_{i=0}^{n}a_i\lambda^{(i+1)\alpha}C(\theta T^*)^{i+1}p_i(t),
$$
and $d = \frac{c^{2\beta}}{2n+2}, a_i=\frac{c^{(i+1)\alpha}}{i+1}$.

Thus \eqref{maineq} is transformed into a Hamiltonian system with impulsive effects:
\begin{equation}\label{XH}
\begin{array}{ll}
 X_{H}:\left\{ \begin{array}{ll}
\theta'=\frac{\partial}{\partial \lambda}H=2\beta d\lambda^{2\beta-1}+\frac{\partial}{\partial \lambda} h_1+\frac{\partial }{\partial \lambda}h_2,\ \ \\[0.2cm]
\lambda'=-\frac{\partial}{\partial \theta}H=-\frac{\partial }{\partial \theta}h_1-\frac{\partial}{\partial \theta} h_2, \quad \ \ \ \ t\neq t_{j};\\[0.2cm]
\lambda(t_{j}+)=\lambda(t_{j}-),\\[0.2cm]
\theta(t_{j}+)=-\theta(t_{j}-),\ \ \ \ \ \ \quad j=\pm1,\pm2,\cdots.
 \end{array}
 \right.
 \end{array}
 \end{equation}

\section{Canonical transformations}

In order to perform some canonical transformations for Hamiltonian system \eqref{XH}, we first introduce a space of functions $\mathcal{F}^k(r)$ which for $r\in \mathbb{R}$ and $k\in \mathbb{N}$ denotes the sets of functions $f(\lambda,\theta,t)$ smooth in $(\lambda,\theta)$, $C^k$ in $t$, and
$$
\sup_{\lambda\ge\lambda_0,(\theta,t)\in \mathbb{S}^1\times \mathbb{S}^1}\lambda^{j-r}\left|D^j_\lambda D^l_\theta f(\lambda,\theta,t)\right| < +\infty.
$$
From the definition, the following properties can be readily verified.
\begin{lemma}\label{property of f}
$$
\begin{array}{ll}
(i)\ \ \mbox{if}\ r_{1}<r_{2}\ \mbox{then}\ \mathcal{F}^k(r_{1})\subset\mathcal{F}^k(r_{2}),\\[0.2cm]
(ii)\ \mbox{if}\ f\in\mathcal{F}^k(r)\ \mbox{then}\ (D_{\lambda})^{j}f\in\mathcal{F}^k(r-j),\\[0.2cm]
(iii)\ \mbox{if}\ f_{1}\in\mathcal{F}^k(r_{1})\ \mbox{and}\ f_{2}\in\mathcal{F}^k(r_{2})\ \mbox{then}\ f_{1}\cdot f_{2}\in\mathcal{F}^k(r_{1}+r_{2}),\\[0.2cm]
(iv)\ \mbox{if}\ f\in\mathcal{F}^k(r)\ \mbox{satisfies}\ |f(\lambda,\cdot)|\geq c\lambda^{r}\ \mbox{for}\ \lambda>\lambda_{0}\ \mbox{then}\ \frac{1}{f}\in\mathcal{F}^k(-r),\\[0.2cm]
(v)\ \mbox{if}\ f\in \mathcal{F}^k(r),\ g\in \mathcal{F}^k(s), \mbox{then}\ f(g)\in\mathcal{F}^k(rs).
\end{array}
$$
\end{lemma}

In this notation, by \eqref{H} we have
$$
h_1(\lambda,\theta,t)\in \mathcal{F}^1\left(\frac{2n + 1}{n + 2}\right),\quad\quad h_2(\lambda,\theta,t) \in \mathcal{F}^0\left(\frac{n + 1}{n + 2}\right).
$$
For $f\in\mathcal{F}^k(r)$ we denote the meanvalue over the $\theta$-variable by $[f]$:
$$
[f](\lambda,t):=\displaystyle\int_{0}^{1}f(\lambda,\theta,t)d\theta.
$$
If $\lambda_{0}>0$, then denote by $\mathcal{A}_{\lambda_{0}}\subset\mathbb{R}^{+}\times \mathbb{S}^{2}$ where $\mathbb{S}^{2}=\mathbb{R}^2/\mathbb{Z}^2$ the annulus
$$
A_{\lambda_{0}}:=\{(\lambda,\theta,t)|\lambda\geq \lambda_{0}\ \mbox{and}\ (\theta,t)\in \mathbb{S}^{2}\}.
$$

For general Hamiltonian functions corresponding to Hamilton system without impulsive effects, there exists a canonical diffeomorphism such that after several times of using this transformation some terms containing angle-variable eventually belong to $\mathcal{F}^k(-\varepsilon)$ for given $\varepsilon$. This work has been well done by R. Dieckerhoff and E. Zehnder in \cite{Dieckerhoff}. Accordingly, for Hamilton system with impulsive effects, we also need a similar canonical transformation to simplify the high order terms containing angle-variable. The vital important point is how to deal with the transformation at impulsive times. Fortunately, due to the special impulsive functions chosen in \eqref{XH}, we find that under the similar canonical transformation as in \cite{Dieckerhoff} the changes of new variables at impulsive times are invariant with the changes of old ones. Such a canonical transformation is what we need. In the following proposition we introduce the canonical transformation which transforms the Hamiltonian $H$ defined by \eqref{H} into new ones, which is more closer to an integrable Hamiltonian.

For this purpose, noticing that $h_1$ $h_2$ given by \eqref{H} are even in angle-variable $\theta$ due to property $(iv)$ of $C(t)$, we consider a slightly general Hamiltonian
\begin{equation}\label{generalH}
H(\lambda,\theta,t) = h_0(\lambda,t) + h_1(\lambda,\theta,t) + h_2(\lambda,\theta,t),\ \ t\neq t_{j},
\end{equation}
with impulsive effects
\begin{equation}\label{impulsive effect old}
\lambda(t_{j}+)=\lambda(t_{j}-),\ \ \ \ \theta(t_{j}+)=-\theta(t_{j}-),\ \ \ \quad j=\pm1,\pm2,\cdots,
\end{equation}
where
$$
h_0\in \mathcal{F}^1\left(\frac{2n+2}{n+2}\right),\quad h_1\in \mathcal{F}^1(a_1),\quad h_2\in \mathcal{F}^0(a_2)
$$
and
$$
D_{\lambda} h_0 \ge c\cdot\lambda^{\frac{2n+2}{n+2}-1}>0,\quad a_1\le\frac{2n+1}{n+2},\quad a_2=\frac{n+1}{n+2}.
$$
In addition, assume that $h_1$, $h_2$ are even in $\theta$.

\begin{proposition}\label{Prop}
There exists a canonical transformation $\psi$ depending periodically on $t$ for $t\neq t_{j},\ j=\pm1,\pm2,\cdots$ of the form
$$
\begin{array}{ll}
\psi:\quad \begin{array}{ll}
 \lambda=\mu+u(\mu,\phi,t),\ \ \ \ \\[0.2cm]
 \theta=\phi+v(\mu,\phi,t),\ \ \  \\[0.2cm]
\end{array}
 \end{array}
$$
with $u\in\mathcal{F}^1\left(a_1-\frac{n}{n+2}\right)$ and $v \in \mathcal{F}^1\left(a_1-\frac{2n+2}{n + 2}\right)$ such that $\mathcal{A}_{\mu^{+}}\subset\psi(\mathcal{A}_{\mu_{0}})\subset\mathcal{A}_{\mu^{-}}$ for some large $\mu^{-}<\mu_{0}<\mu^{+}$.
The transformed Hamiltonian is of the form
\begin{equation}\label{newH}
\widetilde{H}(\mu,\phi,t) = \tilde{h}_0(\mu,t) + \tilde{h}_1(\mu,\phi,t) + \tilde{h}_2(\mu,\phi,t), \ \ \ t\neq t_{j},
\end{equation}
with impulsive effects
$$
\mu(t_{j}+)=\mu(t_{j}-),\ \ \ \ \phi(t_{j}+)=-\phi(t_{j}-),\ \ \ \quad j=\pm1,\pm2,\cdots,
$$
where
$$
\tilde{h}_0\in \mathcal{F}^1\left(\frac{2n + 2}{n + 2}\right),\quad \tilde{h}_1\in \mathcal{F}^1\left(2a_1-\frac{2n+2}{n+2}\right),\quad \tilde{h}_2\in \mathcal{F}^0\left(a_2\right),
$$
$$
\tilde{h}_0(\mu,t)= h_0(\mu,t)+[h_1](\mu,t)\ \ \ \ t\neq t_{j},\ j=\pm1,\pm2,\cdots.
$$
Moreover, the function $\tilde h_0(\mu,t)$ satisfies
$$
D_\mu\tilde h_0\ge c\cdot\mu^{\frac{2n+2}{n+2}-1}>0.
$$
\end{proposition}

\Proof We define the canonical transformation implicitly by
\begin{equation}\label{implicit transform}
\begin{array}{ll}
\psi:\quad \begin{array}{ll}
 \lambda=\mu + D_\theta S(\mu,\theta,t) =:\mu+ \nu,\ \ \ \ \\[0.2cm]
 \phi = \theta + D_\mu S(\mu,\theta,t),\  \\[0.2cm]
\end{array}
 \end{array}
\end{equation}
for $t\neq t_{j},\ j=\pm1,\pm2,\cdots$, where the generating function $S$ will be determined later. The transformation at impulsive times will be calculated by \eqref{impulsive effect old} and $S$ after $S$ is solved. Under this transformation,
for $t\neq t_{j}$ the new Hamiltonian is
$$
\begin{array}{lll}
\widetilde{H}(\mu,\theta,t) &=& h_0(\mu+ \nu, t) + h_1(\mu+ \nu, \theta, t) + h_2(\mu+ \nu, \theta, t)+ D_t S(\mu,\theta,t)\\[0.2cm]
&=& h_0(\mu,t) + D_\mu h_0(\mu,t)\nu + \int^1_0 (1-s)D^2_\mu h_0(\mu + s\nu,t)\nu^2ds\\[0.2cm]
&&+ h_1(\mu, \theta, t)+\int^1_0 D_\mu h_1(\mu + s\nu,\theta, t)\nu ds\\[0.2cm]
&&+ h_2(\mu+ \nu, \theta, t)+ D_t S(\mu,\theta,t).
\end{array}
$$

Let
\begin{equation}\label{mu}
\nu(\mu,\theta,t)= \frac{h_1-[h_1]}{D_\mu h_0}\in\mathcal{F}^1\left(a_1-\frac{n}{n+2}\right),\ \ \ \ t\neq t_{j},
\end{equation}
and
\begin{equation}\label{S}
S(\mu,\theta,t)=\int_0^\theta\,\nu(\mu,s,t)ds,\ \ \ \ t\neq t_{j}.
\end{equation}
Then we obtain for $ t\neq t_{j}$,
$$
\begin{array}{lll}
\widetilde{H}(\mu,\theta,t) &=& h_0(\mu,t) + [h_1](\mu,t)\\[0.2cm]
&&+ \int^1_0 (1-s)D^2_\mu h_0(\mu + s\nu,t)\nu^2ds+\int^1_0 D_\mu h_1(\mu + s\nu,\theta, t)\nu ds\\[0.2cm]
&&+ h_2(\mu+ \nu, \theta, t)+ D_t S(\mu,\theta,t)\\[0.2cm]
&:=&\bar{h}_0(\mu,t)+\bar{h}_1(\mu,\theta, t)+\bar{h}_2(\mu, \theta, t),
\end{array}
$$
where
$$
\begin{array}{lll}
\bar{h}_0(\mu,t) &=& h_0(\mu,t) + [h_1](\mu,t),\\[0.2cm]
\bar{h}_1(\mu,\theta,t) &=& h_0(\mu+ \nu, t)-h_0(\mu,t) - D_\mu h_0(\mu,t)\nu+h_1(\mu+ \nu, \theta, t)-h_1(\mu, \theta, t)\\[0.2cm]
&=&\int^1_0 (1-s)D^2_\mu h_0(\mu + s\nu,t)\nu^2ds+\int^1_0 D_\mu h_1(\mu + s\nu,\theta, t)\nu ds,\\[0.2cm]
\bar{h}_2(\mu,\theta,t) &=&  h_2(\mu+ \nu, \theta, t)+ D_t S(\mu,\theta,t).
\end{array}
$$

It is easy to see that
$$
\bar{h}_0\in \mathcal{F}^1\left(\frac{2n + 2}{n + 2}\right).
$$
 Moreover,
$$
\begin{array}{lll}
\bar{h}_1(\mu,\theta,t) &=& \dint_0^1(1-s)D_\mu^2h_0(\mu+s\nu,t)\nu^2ds + \dint_0^1D_\mu h_1(\mu+s\nu,\theta,t)\nu ds\\[0.42cm]
&\in& \mathcal{F}^1\left(2a_1-\frac{2n+2}{n+2}\right)
\end{array}
$$
and
$$
\bar{h}_2(\mu,\theta,t)\in\mathcal{F}^0(a_2).
$$

From the equation $ \phi= \theta+ D_\mu S(\mu,\theta,t)$, by contraction principle we can solve that $\theta = \phi+v(\mu,\phi,t)$ with $v \in \mathcal{F}^1\left(a_1-\frac{2n+2}{n + 2}\right)$ for large $\mu$, where
\begin{equation}\label{v}
v=-D_\mu S(\mu,\phi+v).
\end{equation}
Then we insert $\theta = \phi+v$ into the first equation of \eqref{implicit transform} and define $u$ to be
\begin{equation}\label{u}
\lambda=\mu+\nu(\mu,\phi+v,t)=\mu+u(\mu,\phi,t).
\end{equation}
Since $\nu\in\mathcal{F}^1\left(a_1-\frac{n}{n+2}\right)$ one concludes that $u\in\mathcal{F}^1\left(a_1-\frac{n}{n+2}\right)$.
Up to now the transformation defined by \eqref{implicit transform} can be rewritten as
$$
\begin{array}{ll}
\psi:\quad \begin{array}{ll}
 \lambda=\mu+u(\mu,\phi,t),\ \ \ \ \\[0.2cm]
 \theta=\phi+v(\mu,\phi,t),\ \ \  \\[0.2cm]
\end{array}
 \end{array}
$$
for $t\neq t_{j},\ j=\pm1,\pm2,\cdots$.
And one verifies easily that the map $\psi$ has a right inverse of the same type as $\psi$ defined on $\mathcal{A}_{\mu^{-}}$ for some large $\mu^{-}$ and that it is injective on $\mathcal{A}_{\mu_{0}}$ for $\mu_{0}$ large. Thus the new Hamiltonian is
$$
\begin{array}{lll}
\widetilde{H}(\mu,\phi,t) &=& \widetilde{H}(\mu,\phi+v(\mu,\phi,t), t)\\[0.2cm]
&=& \bar{h}_0(\mu,t)+ \bar{h}_1(\mu,\phi+v(\mu,\phi,t), t)+\bar{h}_2(\mu, \phi+v(\mu,\phi,t), t)\\[0.2cm]
&:=&\tilde{h}_0(\mu,t) + \tilde{h}_1(\mu,\phi,t) + \tilde{h}_2(\mu,\phi,t),\ \ \ \ t\neq t_{j},
\end{array}
$$
where
$$
\tilde{h}_0\in\mathcal{F}^1\left(\frac{2n + 2}{n + 2}\right),\quad \tilde{h}_1\in \mathcal{F}^1\left(2a_1-\frac{2n+2}{n+2}\right),\quad \tilde{h}_2\in \mathcal{F}^0\left(a_2\right).
$$

For impulsive times, we shall show that under the above transformation $\psi$,
impulsive conditions in \eqref{impulsive effect old} hold true for new variables, that is
$$
\mu(t_{j}+)=\mu(t_{j}-),\ \ \ \phi(t_{j}+)=-\phi(t_{j}-),\ \ \ j=\pm1,\pm2,\cdots.
$$
In fact, from \eqref{implicit transform}, we estimate $\mu(t)$ and $\phi(t)$ at $t=t_{j}$ with
$$
\begin{array}{ll}
 \begin{array}{ll}
 \lambda(t_{j}\pm)=\mu(t_{j}\pm)+D_\theta S(\mu(t_{j}\pm),\theta(t_{j}\pm),t_{j}),\ \ \ \ \\[0.2cm]
 \phi(t_{j}\pm)=\theta(t_{j}\pm)+D_\mu S(\mu(t_{j}\pm),\theta(t_{j}\pm),t_{j}).\ \ \\[0.2cm]
\end{array}
 \end{array}
$$
Combining with the impulsive condition in \eqref{impulsive effect old}
$$
\lambda(t_{j}+)=\lambda(t_{j}-),\ \ \ \
\theta(t_{j}+)=-\theta(t_{j}-),\ \ \ j=\pm1,\pm2,\cdots,
$$
one has
\begin{equation}\label{estimate mu phi}
\begin{array}{ll}
 \begin{array}{ll}
\mu(t_{j}+)+D_ \theta S(\mu(t_{j}+),\theta(t_{j}+),t_{j})=D_\theta S(\mu(t_{j}-),\theta(t_{j}-),t_{j}),\\[0.2cm]
\phi(t_{j}+)-D_\mu S(\mu(t_{j}+),\theta(t_{j}+),t_{j})=-\phi(t_{j}-)+D_\mu S(\mu(t_{j}-),\theta(t_{j}-),t_{j}).
\end{array}
 \end{array}
\end{equation}
Since $h_1$ $h_2$ are even in $\theta$ and by \eqref{mu}
$$
\nu(\mu,\theta,t)= \frac{h_1-[h_1]}{D_\mu h_0},
$$
then $\nu=D_\theta S$ is even and $S=\int_0^\theta\,\nu(\mu,s,t)ds$ is odd in $\theta$. Therefore, the first equation in \eqref{estimate mu phi} becomes
$$
\mu(t_{j}+)+D_\theta S(\mu(t_{j}+),\theta(t_{j}-),t_{j})=\mu(t_{j}-)+D_\theta S(\mu(t_{j}-),\theta(t_{j}-),t_{j}),
$$
which can be written as
$$
\begin{array}{ll}
\big|\mu(t_{j}+)-\mu(t_{j}-)\big|&=\Big|D_\theta S(\mu(t_{j}+),\theta(t_{j}-),t_{j})-D_\theta S(\mu(t_{j}-),\theta(t_{j}-),t_{j})\Big|\\[0.2cm]
&\leq\Big|D_\mu D_ \theta S(\cdot,\theta(t_{j}-),t_{j})\Big|\big|\mu(t_{j}+)-\mu(t_{j}-)\big|.
\end{array}
$$
Because $S\in\mathcal{F}^1\left(a_1-\frac{n}{n+2}\right)$, for large $\mu$, $\big|D_ \mu D_ \theta S\big|$ is small enough such that the above inequality implies
$$
\mu(t_{j}+)=\mu(t_{j}-),\ \ \ j=\pm1,\pm2,\cdots.
$$
Take it back to the second equation in \eqref{estimate mu phi}, then
$$
\begin{array}{ll}
&\phi(t_{j}+)-D_ \mu S(\mu(t_{j}+),-\theta(t_{j}-),t_{j})\\[0.2cm]
=&\phi(t_{j}+)+D_ \mu S(\mu(t_{j}-),\theta(t_{j}-),t_{j})\\[0.2cm]
=&-\phi(t_{j}-)+D_ \mu S(\mu(t_{j}-),\theta(t_{j}-),t_{j}),
\end{array}
$$
from the last two equalities of which one concludes
$$
\phi(t_{j}+)=-\phi(t_{j}-),\ \ \ j=\pm1,\pm2,\cdots.
$$
Thus we have finished the proof of this proposition.\qed

By this Proposition, in new variables $\mu$ and $\phi$, $H$ defined by \eqref{generalH} becomes
$$
\widehat{H}=h_0(\mu+u,t)+h_{1}(\mu+u,\phi+v,t)+h_{2}(\mu+u,\phi+v,t),\ \ \ \ t\neq t_{j}.
$$
We say $\widehat{H}$ is also even in $\phi$ such that this proposition can be applied on $X_{H}$ all the times. Actually, in the proof of Proposition \ref{Prop}, $u\ v$ in transformation $\psi$ are solved by \eqref{u} and \eqref{v}. Since $S$ is odd in $\theta$,
$$
v(-\phi)=-D_\mu S(\mu,-\phi+v(-\phi))= D_\mu S(\mu,\phi-v(-\phi)),
$$
for convenience we suppress variables $\mu$ and $t$ in $v$. Then
$$
\begin{array}{ll}
\big|v(\phi)+v(-\phi)\big|&=\Big|D_\mu S(\mu,\phi+v(\phi))-D_\mu S(\mu,\phi-v(-\phi))\Big|\\[0.2cm]
&\leq\Big|D_\mu D_\theta S\Big|\big|v(\phi)+v(-\phi)\big|,
\end{array}
$$
from which we infer $v$ is odd in $\phi$ for $\mu$ large enough. And the fact $u$ is even in $\phi$ is easily obtained by $u=\nu(\mu,\phi+v)$ where $\nu$ is even and $v$ is odd. Therefore, combining the odevity of $u$, $v$ $h_{1}$ and, $h_{2}$,
$$
\begin{array}{ll}
&\widehat{H}(-\phi)\\[0.2cm]
=&h_0(\mu+u(-\phi),t)+h_{1}(\mu+u(-\phi),-\phi+v(-\phi),t)+h_{2}(\mu+u(-\phi),-\phi+v(-\phi),t)\\[0.2cm]
=&\widehat{H}(\phi).
\end{array}
$$
Applying this proposition many times to the original system (\ref{H}), the transformed Hamiltonian is of the form
\begin{equation}\label{H after several tran}
H(\lambda,\theta,t) = h_0(\lambda,t) + h_1(\lambda,\theta,t) + h_2(\lambda,\theta,t),\ \ \ \ t\neq t_{j},
\end{equation}
with impulsive effects
$$
\lambda(t_j+)=\lambda(t_j-),\ \ \ \theta(t_j+)=-\theta(t_j-),\ \ \ \ j=\pm1,\pm2,\cdots,
$$
which has the property
$$
h_0\in \mathcal{F}^1\left(\frac{2n + 2}{n + 2}\right),\quad h_1\in \mathcal{F}^1\left(\frac{1-4n}{n + 2}\right),\quad h_2\in \mathcal{F}^0\left(\frac{n + 1}{n + 2}\right).
$$

Because $h_2\in \mathcal{F}^0\left(\frac{n + 1}{n + 2}\right)$, let $\bar{h}_2(\lambda,\theta,t) = \lambda^{-\frac{n+1}{n+2}}h_2(\lambda,\theta,t)$, then there is a constant $M>0$ such that, for $\lambda\geq \lambda_0 ,(\theta,t)\in \mathbb{S}^1\times \mathbb{S}^1$
$$
|\lambda^jD^j_\lambda D^l_\theta\bar{h}_2(\lambda,\theta,t)|\le M,
$$
where $j+l$ is bounded by a sufficiently large integer $K>0$, which guarantees the smoothness assumption $C^4$ in the large twist theorem for  the Poincar\'{e} map, hence we can choose the constant $M$ independently of $j,l$ for $j+l\leq K$.

For $j+l\leq K$, we expand the function $\bar{h}_{2}^{jl}(\lambda,\theta,t):=\lambda^jD^j_\lambda D^l_\theta\bar{h}_2(\lambda,\theta,t)$ into the Fourier series with respect to $t$
$$
\frac{\bar{h}_{2}^{jl0}}{2}+\sum_{n=1}^{+\infty}\left(\bar{h}_{2}^{jlnc}\cos 2\pi nt+\bar{h}_{2}^{jlns}\sin 2\pi nt\right)
$$
with the Fourier coefficients $\bar{h}_{2}^{jlnc}, \bar{h}_{2}^{jlns}$,  and its part sum is
$$
S_n(\bar{h}_{2}^{jl})=\frac{\bar{h}_{2}^{jl0}}{2}+\sum_{k=1}^{n}\left(\bar{h}_{2}^{jlkc}\cos 2\pi kt+\bar{h}_{2}^{jlks}\sin 2\pi kt\right).
$$

Consider the Fej\'{e}r sum of $\bar{h}_{2}^{jl}$
$$
F_n(\bar{h}_{2}^{jl})=\frac{\dsum_{k=0}^n S_k(\bar{h}_{2}^{jl})}{n+1}.
$$
By Fej\'{e}r Theorem,  the Fej\'{e}r sum $F_n(\bar{h}_{2}^{jl})$ of $\bar{h}_{2}^{jl}$ converges to $\bar{h}_{2}^{jl}$ uniformly with respect to $t$.
For any fixed $\epsilon > 0$, there exists a sufficiently large $N_{jl}\in \mathbb{N}$ such that for $n\geq N$ and $\lambda\geq \lambda_0 ,(\theta,t)\in \mathbb{S}^1\times \mathbb{S}^1$,
$$
|\bar{h}_{2}^{jl}(\lambda,\theta,t)-F_n(\bar{h}_{2}^{jl})(\lambda,\theta,t)|<\epsilon.
$$
Since the integers $j,l$ satisfying $j+l\leq K$ are finite, the positive integer $N_{jl}$ can be chosen independently of $j, l$. Let
$$
\bar{h}_{21}(\lambda,\theta,t)=F_N(\bar{h}_{2}^{00})(\lambda,\theta,t), \ \ \ \ \bar{h}_{22}(\lambda,\theta,t)=\bar{h}_{2}^{00}(\lambda,\theta,t)-F_N(\bar{h}_{2}^{00})(\lambda,\theta,t).
$$
Then
$$
\bar{h}_2(\lambda,\theta,t) = \bar{h}_{21}(\lambda,\theta,t) + \bar{h}_{22}(\lambda,\theta,t),
$$
and the functions $\bar{h}_{21}$ and $\bar{h}_{22}$ satisfy
 for $\lambda\geq \lambda_0 ,(\theta,t)\in \mathbb{S}^1\times \mathbb{S}^1$, and $j+1\leq K$
$$
|\lambda^jD^j_\lambda D^l_\theta\bar{h}_{21}(\lambda,\theta,t)|\le M+1,
$$
$$
|\lambda^jD^j_\lambda D^l_\theta\bar{h}_{22}(\lambda,\theta,t)|<\epsilon.
$$

From the above discussions, we may rewrite the Hamiltonian (\ref{H after several tran}) in the form
\begin{equation}\label{H3}
H(\lambda,\theta,t) = h_0(\lambda,t) + h_1(\lambda,\theta,t) + h_{21}(\lambda,\theta,t) + h_{22}(\lambda,\theta,t),\ \ \ \ t\neq t_{j}
\end{equation}
where
$$
h_0\in \mathcal{F}^1\left(\frac{2n+2}{n+2}\right),\quad h_1\in \mathcal{F}^1\left(\frac{1-4n}{n+2}\right),
$$
$$
h_{21}\in \mathcal{F}^\infty\left(\frac{n+1}{n+2}\right), \quad h_{22}\in \mathcal{F}^0\left(\frac{n+1}{n+2}\right),
$$
and
$$
\sup_{\lambda\ge\lambda_0,(\theta,t)\in \mathbb{S}^1\times \mathbb{S}^1}\left|\lambda^{j-\frac{n+1}{n+2}}D^j_\lambda D^l_\theta h_{22}(\lambda,\theta,t)\right|<\epsilon.
$$

Applying Proposition \ref{Prop} to the Hamiltonian (\ref{H3}), we may assume that the transformed hamiltonian of (\ref{H3}) is of the form
\begin{equation}\label{H4}
H(\lambda,\theta,t) = h_0(\lambda,t) + h_1(\lambda,\theta,t) + h_2(\lambda,\theta,t),\ \ \ \ t\neq t_{j},
\end{equation}
where the functions $h_0,\ h_1$ and $h_2$ satisfy
$$
h_0(\lambda,t) = d\cdot\lambda^{\frac{2n+2}{n+2}} + \bar h_0(\lambda, t),
$$
and
$$
\bar h_0\in \mathcal{F}^1\left(\frac{2n+1}{n+2}\right),\quad h_1\in \mathcal{F}^1\left(\frac{1-4n}{n+2}\right),\quad h_2\in \mathcal{F}^0\left(\frac{n+1}{n+2}\right).
$$
Moreover,
$$
\sup_{\lambda\ge\lambda_0,(\theta,t)\in \mathbb{S}^1\times \mathbb{S}^1}\left|\lambda^{j-\frac{n+1}{n+2}}D^j_\lambda D^l_\theta h_{2}(\lambda,\theta,t)\right|<\epsilon.
$$

\section{Estimation of the time 1 map of Hamiltonian system with impulsive effects}

Due to the discussion in Section 3, now let $H$ be as in \eqref{H4}, that is
$$
H(\lambda,\theta,t) = h_0(\lambda,t) + h_1(\lambda,\theta,t) + h_2(\lambda,\theta,t),\ \ \ \ t\neq t_{j},
$$
with impulsive effects
$$
\lambda(t_{j}+)=\lambda(t_{j}-),\ \ \ \
\theta(t_{j}+)=-\theta(t_{j}-),\ \ \ j=\pm1,\pm2,\cdots,
$$
where the functions $h_0,\ h_1$ and $h_2$ satisfy
$$
h_0(\lambda,t) = d\cdot\lambda^{\frac{2n+2}{n+2}} + \bar h_0(\lambda, t),
$$
and
$$
\bar h_0\in \mathcal{F}^1\left(\frac{2n+1}{n+2}\right),\quad h_1\in \mathcal{F}^1\left(\frac{1-4n}{n+2}\right),\quad h_2\in \mathcal{F}^0\left(\frac{n+1}{n+2}\right).
$$
Moreover,
$$
\sup_{\lambda\ge\lambda_0,(\theta,t)\in \mathbb{S}^1\times \mathbb{S}^1}\left|\lambda^{j-\frac{n+1}{n+2}}D^j_\lambda D^l_\theta h_{2}(\lambda,\theta,t)\right|<\epsilon.
$$
The corresponding Hamiltonian system is
\begin{equation}\label{XH4}
\begin{array}{ll}
\left\{ \begin{array}{ll}
\dot{\lambda}=-D_{\theta}h_1(\lambda,\theta,t)-D_{\theta}h_{2}(\lambda,\theta,t),\\[0.2cm]
\dot{\theta}=D_{\lambda}h_0(\lambda,t)+D_{\lambda}h_1(\lambda,\theta,t)+D_{\lambda}h_{2}(\lambda,\theta,t),\ \ \ \ t\neq t_{j};\\[0.2cm]
\lambda(t_{j}+)=\lambda(t_{j}-),\\[0.2cm]
\theta(t_{j}+)=-\theta(t_{j}-),\quad \ \ \ \ j=\pm1,\pm2,\cdots.\\[0.2cm]
 \end{array}
\right.
 \end{array}
 \end{equation}
Define
\begin{equation}\label{rho lambda}
\rho=\lambda^{\frac{1}{n+2}}, \quad\quad \theta=\theta,
\end{equation}
then for $t\neq t_{j}$,
$$
\begin{array}{lll}
\frac{d\rho}{dt}&=&\frac{1}{n+2}\lambda^{-\frac{n+1}{n+2}}\frac{d\lambda}{dt}\\[0.2cm]
&=&-\frac{1}{n+2}\lambda^{-\frac{n+1}{n+2}}D_{\theta}h_1(\lambda,\theta,t)-\frac{1}{n+2}\lambda^{-\frac{n+1}{n+2}}D_{\theta}h_{2}(\lambda,\theta, t)\\[0.2cm]
&\triangleq&f_{1}(\rho,\theta,t)+f_{2}(\rho,\theta, t),\\[0.2cm]
\frac{d\theta}{dt}&=&D_{\lambda}h_0(\lambda,t)+D_{\lambda}h_1(\lambda,\theta,t)+D_{\lambda}h_{2}(\lambda,\theta, t)\\[0.2cm]
&\triangleq&f_{3}(\rho,t)+f_{4}(\rho,\theta,t)+f_{5}(\rho,\theta, t),
\end{array}
$$
that is,
\begin{equation*}
\left\{\begin{array}{lll}
\dot{\rho}=f_{1}(\rho,\theta,t)+f_{2}(\rho,\theta, t),\\[0.2cm]
\dot{\theta}=f_{3}(\rho,t)+f_{4}(\rho,\theta,t)+f_{5}(\rho,\theta, t),\ \ \ \ t\neq t_{j};\\[0.2cm]
\rho(t_{j}+)=\rho(t_{j}-),\\[0.2cm]
\theta(t_{j}+)=-\theta(t_{j}-),\quad \ \ \ \ j=\pm1,\pm2,\cdots,
 \end{array}\right.
\end{equation*}
where
$$
f_1\in \mathcal{F}^1(-5n), f_2\in \mathcal{F}^0(0),f_3\in \mathcal{F}^1(n), f_4\in \mathcal{F}^1(-5n-1), f_5\in \mathcal{F}^0(-1).
$$
Moreover,
$$
f_3 = d\cdot\rho^n + O(\rho^{n-1}),\quad |\rho^jD^j_\rho D^l_\theta f_2|<\epsilon,\ \  |\rho^{j+1}D^j_\rho D^l_\theta f_5| <\epsilon.
$$

Now we replace $\rho$ by
$$
\gamma^{-1}I,\quad\quad I\in[1,2],
$$
then $\gamma$ is a parameter with $\rho\to+\infty \Leftrightarrow \gamma\to 0^+$, and the corresponding system has the form
\begin{equation}\label{I theta}
\begin{array}{lll}
\dot{I}&=&\gamma f_{1}(\gamma^{-1}I,\theta,t)+\gamma f_{2}(\gamma^{-1}I,\theta, t),\\[0.2cm]
       &\triangleq&\gamma^{5n+1} \bar{f}_{1}(I,\theta,t;\gamma)+\gamma \bar{f}_{2}(I,\theta,t;\gamma),\\[0.2cm]
\dot{\theta}&=&f_{3}(\gamma^{-1}I,t)+f_{4}(\gamma^{-1}I,\theta,t)+f_{5}(\gamma^{-1}I,\theta, t)\\[0.2cm]
&\triangleq&\gamma^{-n}\bar{f}_{3}(I,t;\gamma)+\gamma^{5n+1}\bar{f}_{4}(I,\theta,t;\gamma)+\gamma \bar{f}_{5}(I,\theta, t;\gamma),\ \ \ \ t\neq t_{j},
 \end{array}
 \end{equation}
and
$$
I(t_{j}+)=I(t_{j}-),\ \ \
\theta(t_{j}+)=-\theta(t_{j}-),\quad \ \ \ \ j=\pm1,\pm2,\cdots,
$$
where $\bar{f}_i$ $(i=1,2,3,4,5)$ satisfy
\begin{equation}\label{A0}
|D_I^k\bar{f}_3(I,t;\gamma)|\leq C,
\end{equation}
\begin{equation}\label{A1}
|D_I^k D_{\theta}^j \bar{f}_1(I,\theta,t;\gamma)|\leq C,\ \  |D_I^k D_{\theta}^j \bar{f}_4(I,\theta,t;\gamma)|\leq C
\end{equation}
\begin{equation}\label{A2}
|D_I^k D_{\theta}^j \bar{f}_2(I,\theta,t;\gamma)| < \epsilon,\ \  |D_I^k D_{\theta}^j \bar{f}_5(I,\theta,t;\gamma)| < \epsilon.
\end{equation}
for some constant $C>0$ and $k+j\leq 4$.

Since the estimations of $\bar{f}_i$ $(i=1,2,3,4,5)$, one verifies easily that the solutions of \eqref{I theta} do exist for $t\in[0,1]$ with $\gamma$ sufficiently small. Next we consider system \eqref{I theta} on $[0,1]$, which includes two impulsive times $t_{1}$ and $t_{2}$ with $t_{2}-t_{1}\neq\frac{1}{2}$. Let $(I(t),\theta(t))$ be the solution of (\ref{I theta}) with the initial condition $(I(0),\theta(0))=(I_{0}, \theta_{0})$, and denote
\begin{equation}
\begin{array}{ll}
\ \begin{array}{ll}
P_{0}:(I_{0},\theta_{0})\mapsto(I(t_{1}-),\theta(t_{1}-)),\\[0.2cm]
P_{1}:(I(t_{1}+),\theta(t_{1}+))\mapsto(I(t_{2}-),\theta(t_{2}-)),\\[0.2cm]
P_{2}:(I(t_{2}+),\theta(t_{2}+))\mapsto(I(1),\theta(1)),\\[0.2cm]
J_{i}:(I(t_{i}-),\theta(t_{i}-))\mapsto(I(t_{i}+),\theta(t_{i}+)),\ i=1,2.
 \end{array}
 \end{array}
 \end{equation}
Then the Poincar\'{e} map of \eqref{I theta} $P:(I_{0},\theta_{0})\mapsto(I(1),\theta(1))$ can be expressed by
$$
P=P_{2}\circ J_{2}\circ P_{1}\circ J_{1}\circ P_{0}.
$$
Since $P_{j},\ j=0,1,2$ are symplectic by the equation $x''+x^{2n+1}+\sum_{j=0}^{2n}x^{j}p_{j}(t)=0$ being conservative and the transformations $\varphi$ and $\psi$ being symplectic, together with $J_{i},\ i=1,2$ being area-preserving homeomorphisms by
$$
\Big|\frac{\partial(I(t_{i}+),\theta(t_{i}+))}{\partial(I(t_{i}-),\theta(t_{i}-))}\Big|=1,
$$
the Poincar\'{e} map $P$ is an area-preserving homeomorphism. Then we estimate $P$ in the next Lemma.

\begin{lemma}\label{estimation}
The time 1 map $P$ of the flow $\phi^{t}$ of the vectorfield given by \eqref{I theta} is of form:
\begin{equation}\label{estimation of P}
\begin{array}{ll}
P:\quad \begin{array}{ll}
I(1)=I_{0}+\Xi_1(I_0,\theta_0; \gamma),\ \ \ \ \\[0.2cm]
\theta(1)=\theta_{0} +\gamma^{-n}\Pi(I_{0};\gamma)+ \Xi_2(I_0,\theta_0; \gamma),\ \ \\[0.2cm]
\end{array}
 \end{array}
 \end{equation}
where
$$
\Pi'(I_{0};\gamma)\geq \frac{d}{2}[1-2(t_2-t_1)],
$$
if $[1-2(t_2-t_1)]>0$ or
$$
\Pi'(I_{0};\gamma)\leq \frac{d}{2}[1-2(t_2-t_1)],
$$
if $[1-2(t_2-t_1)]<0$ for $I\in [1,2]$ with $\gamma$ sufficiently small, and
\begin{equation}\label{Small1}
|D_{I_{0}}^k D_{\theta_{0}}^l\Xi_i(I_0,\theta_0; \gamma)|\leq O(\gamma^{n+1})+\gamma^{-4n}\epsilon \,O(\gamma)
\end{equation}
as $\gamma\to 0^+$ for $i=1,2$ and $k+l\leq 4$.
\end{lemma}
\Proof
We estimate $(I(t),\theta(t))$ on $[0,t_{1}),\ [t_{1},t_{2})$ and $[t_{2},1]$ respectively.
First, for $t\in[0,t_{1})$,
set for the flow $(I(t),\theta(t))=\phi^{t}(I_0,\theta_0)$ with $\phi^{0}=id$
$$
\begin{array}{lll}
I(t)&=&I_0+A_1(I_0,\theta_0,t;\gamma),\\[0.2cm]
\theta(t)&=&\theta_0+\gamma^{-n} B_{11}(I_0,\theta_0,t;\gamma)+B_{12}(I_0,\theta_0,t;\gamma),
\end{array}
$$
where
$$
\begin{array}{lll}
B_{11}(I_0,\theta_0,t;\gamma)&=&\int_0^t\bar{f}_{3}(I(s),s;\gamma)ds,\\[0.2cm]
A_1(I_0,\theta_0,t;\gamma)&=&\gamma^{5n+1} \int_0^t\bar{f}_{1}(I(s),\theta(s),s;\gamma)ds+\gamma \int_0^t\bar{f}_{2}(I(s),\theta(s),s;\gamma)ds,\\[0.2cm]
B_{12}(I_0,\theta_0,t;\gamma)&=&\gamma^{5n+1} \int_0^t\bar{f}_{4}(I(s),\theta(s),s;\gamma)ds+\gamma \int_0^t\bar{f}_{5}(I(s),\theta(s),s;\gamma)ds.
\end{array}
$$
Then by the integral equation of \eqref{I theta}, functions $A_1,B_{11},B_{12}$ satisfy
\begin{equation}\label{ABC1}
\begin{array}{lll}
B_{11}&=&\int_0^t\bar{f}_{3}(I_0+A_1,s;\gamma)ds\\[0.2cm]
A_1&=&\gamma^{5n+1} \int_0^t\bar{f}_{1}(I_0+A_1,\theta_0+\gamma^{-n}B_{11}+B_{12},s;\gamma)ds\\[0.2cm]
&&+\gamma \int_0^t\bar{f}_{2}(I_0+A_1,\theta_0+\gamma^{-n}B_{11}+B_{12},s;\gamma)ds,\\[0.2cm]
B_{12}&=&\gamma^{5n+1} \int_0^t\bar{f}_{4}(I_0+A_1,\theta_0+\gamma^{-n}B_{11}+B_{12},s;\gamma)ds\\[0.2cm]
&&+\gamma \int_0^t\bar{f}_{5}(I_0+A_1,\theta_0+\gamma^{-n}B_{11}+B_{12},s;\gamma)ds.
\end{array}
\end{equation}
Therefore
\begin{equation}\label{Pt1left}
\begin{array}{ll}
\quad \begin{array}{ll}
I(t_1-)=I_{0}+\Xi_{11}(I_0,\theta_0; \gamma),\ \ \ \ \\[0.2cm]
\theta(t_1-)=\theta_{0} +\gamma^{-n}\psi_1(I_{0};\gamma)+ \Xi_{12}(I_0,\theta_0; \gamma),\ \ \\[0.2cm]
\end{array}
 \end{array}
 \end{equation}
where
\begin{equation}\label{psi1}
\begin{array}{lll}
\psi_1(I_0;\gamma)&=&\int_0^{t_1}\,\bar{f}_{3}(I(t),t;\gamma)dt\\[0.2cm]
&=&\int_0^{t_1}\,[d \cdot I(t)^n+ I(t)^{n-1}O(\gamma)]dt\\[0.2cm]
&=&t_1d \cdot (I_0+O(\gamma^{5n+1})+\epsilon O(\gamma))^n+t_1(I_0+O(\gamma^{5n+1})+\epsilon O(\gamma))^{n-1}O(\gamma )\\[0.2cm]
&=&t_1d \cdot I_0^n+t_1O(\gamma)(1+\epsilon)I_0^{n-1}+\cdots,\\[0.2cm]
\end{array}
\end{equation}
and
$$
\begin{array}{lll}
 \Xi_{11}(I_0,\theta_0; \gamma)&=&A_1(I_0,\theta_0,t_1;\gamma)\\[0.2cm]
 &=&\gamma^{5n+1} \int_0^{t_1}\bar{f}_{1}(I_0+A_1,\theta_0+\gamma^{-n}B_{11}+B_{12},s;\gamma)ds\\[0.2cm]
 &&+\gamma \int_0^{t_1}\bar{f}_{2}(I_0+A_1,\theta_0+\gamma^{-n}B_{11}+B_{12},s;\gamma)ds,\\[0.2cm]
 \Xi_{12}(\theta_0, I_0; \gamma)&=&B_{12}(I_0,\theta_0,t_1;\gamma)\\[0.2cm]
 &=&\gamma^{5n+1} \int_0^{t_1}\bar{f}_{4}(I_0+A_1,\theta_0+\gamma^{-n}B_{11}+B_{12},s;\gamma)ds\\[0.2cm]
 &&+\gamma \int_0^{t_1}\bar{f}_{5}(I_0+A_1,\theta_0+\gamma^{-n}B_{11}+B_{12},s;\gamma)ds.
\end{array}
$$
According to (\ref{A1}) and (\ref{A2}), it is easy to see that
\begin{equation}\label{xi1}
|\Xi_{1i}(I_0,\theta_0; \gamma)|\leq O(\gamma^{5n+1})+\epsilon \,O(\gamma)
\end{equation}
as $\gamma\to 0^+$ for $i=1,2$. As for the estimates on the derivatives, from the equations (\ref{ABC1}), we find that
\begin{equation*}
\begin{array}{lll}
D_{I_0}B_{11}&=&\int_0^tD_{I_0}\bar{f}_{3}\cdot(1+D_{I_0}A_1)ds\\[0.2cm]
D_{I_0}A_1&=&\gamma^{5n+1} \int_0^t\left[D_{I_0}\bar{f}_{1}\cdot(1+D_{I_0}A_1)+D_{\theta_0}\bar{f}_{1}\cdot(\gamma^{-n}D_{I_0}B_{11}+D_{I_0}B_{12})\right]ds\\[0.2cm]
&&+\gamma \int_0^t\left[D_{I_0}\bar{f}_{2}\cdot(1+D_{I_0}A_1)+D_{\theta_0}\bar{f}_{2}\cdot(\gamma^{-n}D_{I_0}B_{11}+D_{I_0}B_{12})\right]ds,\\[0.2cm]
D_{I_0}B_{12}&=&\gamma^{5n+1} \int_0^t\left[D_{I_0}\bar{f}_{4}\cdot(1+D_{I_0}A_1)+D_{\theta_0}\bar{f}_{4}\cdot(\gamma^{-n}D_{I_0}B_{11}+D_{I_0}B_{12})\right]ds\\[0.2cm]
&&+\gamma \int_0^t\left[D_{I_0}\bar{f}_{5}\cdot(1+D_{I_0}A_1)+D_{\theta_0}\bar{f}_{5}\cdot(\gamma^{-n}D_{I_0}B_{11}+D_{I_0}B_{12})\right]ds,
\end{array}
\end{equation*}
which together with (\ref{A0}), (\ref{A1}), (\ref{A2}) implies that
$$|D_{I_0}\Xi_{1i}(I_0,\theta_0; \gamma)|\leq O(\gamma^{4n+1})+\gamma^{-n}\epsilon \,O(\gamma)$$
for $i=1,2$. Differentiating $A_1,B_{11},B_{12}$ with respect to $I_0$ $k$ times yields that
$$|D_{I_0}^k\Xi_{1i}(I_0,\theta_0; \gamma)|\leq O(\gamma^{5n+1-kn})+\gamma^{-kn}\epsilon \,O(\gamma)$$
for $i=1,2$ and $0\leq k\leq 4$. Also by Differentiating $A_1,B_{11},B_{12}$ with respect to $\theta_0$, we know that
\begin{equation*}
\begin{array}{lll}
D_{\theta_0}B_{11}&=&\int_0^tD_{I_0}\bar{f}_{3}\cdot D_{\theta_0}A_1ds\\[0.2cm]
D_{\theta_0}A_1&=&\gamma^{5n+1} \int_0^t\left[D_{I_0}\bar{f}_{1}\cdot D_{\theta_0}A_1+D_{\theta_0}\bar{f}_{1}\cdot(1+\gamma^{-n}D_{\theta_0}B_{11}+D_{\theta_0}B_{12})\right]ds\\[0.2cm]
&&+\gamma \int_0^t\left[D_{I_0}\bar{f}_{2}\cdot D_{\theta_0}A_1+D_{\theta_0}\bar{f}_{2}\cdot(1+\gamma^{-n}D_{\theta_0}B_{11}+D_{\theta_0}B_{12})\right]ds,\\[0.2cm]
D_{\theta_0}B_{12}&=&\gamma^{5n+1} \int_0^t\left[D_{I_0}\bar{f}_{4}\cdot D_{\theta_0}A_1+D_{\theta_0}\bar{f}_{4}\cdot(1+\gamma^{-n}D_{\theta_0}B_{11}+D_{\theta_0}B_{12})\right]ds\\[0.2cm]
&&+\gamma \int_0^t\left[D_{I_0}\bar{f}_{5}\cdot D_{\theta_0}A_1+D_{\theta_0}\bar{f}_{5}\cdot(1+\gamma^{-n}D_{\theta_0}B_{11}+D_{\theta_0}B_{12})\right]ds,
\end{array}
\end{equation*}
which together with (\ref{A0}), (\ref{A1}), (\ref{A2}) implies that
$$|D_{\theta_0}\Xi_{1i}(I_0,\theta_0; \gamma)|\leq O(\gamma^{5n+1})+\epsilon \,O(\gamma), \ \ i=1,2.$$
The estimates on the higher derivatives $D_{I_0}^k D_{\theta_0}^l\Xi_{1i}$ can be proven similarly, from which we have
\begin{equation}\label{xi1}
|D_{I_{0}}^k D_{\theta_{0}}^l\Xi_{1i}(I_0,\theta_0; \gamma)|\leq O(\gamma^{n+1})+\gamma^{-4n}\epsilon \,O(\gamma)
\end{equation}
as $\gamma\to 0^+$ for $i=1,2$ and $k+l\leq 4$.

By the impulsive functions in \eqref{I theta} and the value at $t_1-$ in \eqref{Pt1left}, there are
\begin{equation}\label{Pt1right}
\begin{array}{ll}
\quad \begin{array}{ll}
I(t_1+)=I(t_1-)=I_{0}+\Xi_{11}(I_0,\theta_0; \gamma),\ \ \ \ \\[0.2cm]
\theta(t_1+)=-\theta(t_1-)=-\theta_{0} -\gamma^{-n}\psi_1(I_{0};\gamma)- \Xi_{12}(I_0,\theta_0; \gamma).\ \ \\[0.2cm]
\end{array}
 \end{array}
 \end{equation}
Then for $t\in[t_1,t_2)$, we assume
$$
\begin{array}{lll}
I(t)&=&I(t_1+)+A_2(I_0,\theta_0,t;\gamma),\\[0.2cm]
\theta(t)&=&\theta(t_1+)+\gamma^{-n} B_{21}(I_0,\theta_0,t;\gamma)+B_{22}(I_0,\theta_0,t;\gamma),
\end{array}
$$
where
$$
\begin{array}{lll}
B_{21}(I_0,\theta_0,t;\gamma)&=&\int_{t_1}^t\bar{f}_{3}(I(s),s;\gamma)ds,\\[0.2cm]
A_2(I_0,\theta_0,t;\gamma)&=&\gamma^{5n+1} \int_{t_1}^t\bar{f}_{1}(I(s),\theta(s),s;\gamma)ds+\gamma \int_{t_1}^t\bar{f}_{2}(I(s),\theta(s),s;\gamma)ds,\\[0.2cm]
B_{22}(I_0,\theta_0,t;\gamma)&=&\gamma^{5n+1} \int_{t_1}^t\bar{f}_{4}(I(s),\theta(s),s;\gamma)ds+\gamma \int_{t_1}^t\bar{f}_{5}(I(s),\theta(s),s;\gamma)ds.
\end{array}
$$
By the integral equation of \eqref{I theta}, functions $A_2,B_{21},B_{22}$ satisfy
\begin{equation}\label{ABC2}
\begin{array}{lll}
B_{21}&=&\int_{t_1}^t\bar{f}_{3}(I({t_1}+)+A_2,s;\gamma)ds\\[0.2cm]
A_2&=&\gamma^{5n+1} \int_{t_1}^t\bar{f}_{1}(I({t_1}+)+A_2,\theta({t_1}+)+\gamma^{-n}B_{21}+B_{22},s;\gamma)ds\\[0.2cm]
&&+\gamma \int_{t_1}^t\bar{f}_{2}(I({t_1}+)+A_2,\theta({t_1}+)+\gamma^{-n}B_{21}+B_{22},s;\gamma)ds,\\[0.2cm]
B_{22}&=&\gamma^{5n+1} \int_{t_1}^t\bar{f}_{4}(I({t_1}+)+A_2,\theta({t_1}+)+\gamma^{-n}B_{21}+B_{22},s;\gamma)ds\\[0.2cm]
&&+\gamma \int_{t_1}^t\bar{f}_{5}(I({t_1}+)+A_2,\theta({t_1}+)+\gamma^{-n}B_{21}+B_{22},s;\gamma)ds.
\end{array}
\end{equation}
Similarly, we have
\begin{equation}\label{Pt2left}
\begin{array}{ll}
\quad \begin{array}{ll}
I(t_2-)&=I(t_1+)+\Xi_{21}(I_0,\theta_0; \gamma)\\[0.2cm]
&=I_{0}+\Xi_{11}(\theta_{0}, I_{0}; \gamma)+\Xi_{21}(I_0,\theta_0; \gamma),\ \ \ \ \\[0.2cm]
\theta(t_2-)&=\theta(t_1+) +\gamma^{-n}\psi_2(I_{0};\gamma)+ \Xi_{22}(I_0,\theta_0; \gamma)\\[0.2cm]
&=-\theta_{0} -\gamma^{-n}\psi_1(I_{0};\gamma)- \Xi_{12}(I_0,\theta_0; \gamma)+\gamma^{-n}\psi_2(I_{0};\gamma)+ \Xi_{22}(I_0,\theta_0; \gamma),\ \ \\[0.2cm]
\end{array}
 \end{array}
 \end{equation}
where
\begin{equation}\label{psi2}
\begin{array}{lll}
\psi_2(I_0;\gamma)&=&\int_{t_1}^{t_2}\,\bar{f}_{3}(I(t),t;\gamma)dt\\[0.2cm]
&=&\int_{t_1}^{t_2}\,[d \cdot I(t)^n+ I(t)^{n-1}O(\gamma)]dt\\[0.2cm]
&=&(t_2-t_1)d \cdot \big(I(t_1+)+O(\gamma^{5n+1})+\epsilon O(\gamma)\big)^n\\[0.2cm]
&&+(t_2-t_1)\big(I(t_1+)+O(\gamma^{5n+1})+\epsilon O(\gamma)\big)^{n-1}O(\gamma )\\[0.2cm]
&=&(t_2-t_1)d \cdot I(t_1+)^n+(t_2-t_1)O(\gamma)(1+\epsilon)I(t_1+)^{n-1}+\cdots\\[0.2cm]
&=&(t_2-t_1)d \cdot \big(I_0+\Xi_{11}\big)^n+(t_2-t_1)O(\gamma)(1+\epsilon)\big(I_0+\Xi_{11}\big)^{n-1}+\cdots\\[0.2cm]
&=&(t_2-t_1)d \cdot \big(I_0+O(\gamma^{5n+1})+\epsilon O(\gamma)\big)^n\\[0.2cm]
&&+(t_2-t_1)O(\gamma)(1+\epsilon)\big(I_0+O(\gamma^{5n+1})+\epsilon O(\gamma)\big)^{n-1}+\cdots\\[0.2cm]
&=&(t_2-t_1)\big(d \cdot I_0^n+O(\gamma)(1+\epsilon)I_0^{n-1}+\cdots\big)
\end{array}
\end{equation}
and
$$
\begin{array}{lll}
 \Xi_{21}(I_0,\theta_0; \gamma)&=&A_2(I_0,\theta_0,t_2;\gamma)\\[0.2cm]
 &=&\gamma^{5n+1} \int_{t_1}^{t_2}\bar{f}_{1}(I(t_1+)+A_2,\theta(t_1+)+\gamma^{-n}B_{21}+B_{22},s;\gamma)ds\\[0.2cm]
 &&+\gamma \int_{t_1}^{t_2}\bar{f}_{2}(I(t_1+)+A_2,\theta(t_1+)+\gamma^{-n}B_{21}+B_{22},s;\gamma)ds\\[0.2cm]
&=&\gamma^{5n+1} \int_{t_1}^{t_2}\bar{f}_{1}(I_0+\Xi_{11}+A_2,-\theta_{0} -\gamma^{-n}\psi_1-\Xi_{12}+\gamma^{-n}B_{21}+B_{22},s;\gamma)ds\\[0.2cm]
 &&+\gamma \int_{t_1}^{t_2}\bar{f}_{2}(I_0+\Xi_{11}+A_2,-\theta_{0} -\gamma^{-n}\psi_1-\Xi_{12}+
 \gamma^{-n}B_{21}+B_{22},s;\gamma)ds,\\[0.2cm]
 \Xi_{22}(I_0,\theta_0; \gamma)&=&B_{22}(I_0,\theta_0,t_2;\gamma)\\[0.2cm]
 &=&\gamma^{5n+1} \int_{t_1}^{t_2}\bar{f}_{4}(I(t_1+)+A_2,\theta(t_1+)+\gamma^{-n}B_{21}+B_{22},s;\gamma)ds\\[0.2cm]
 &&+\gamma \int_{t_1}^{t_2}\bar{f}_{5}(I(t_1+)+A_2,\theta(t_1+)+\gamma^{-n}B_{21}+B_{22},s;\gamma)ds\\[0.2cm]
&=&\gamma^{5n+1} \int_{t_1}^{t_2}\bar{f}_{4}(I_0+\Xi_{11}+A_2,-\theta_{0} -\gamma^{-n}\psi_1-\Xi_{12}+\gamma^{-n}B_{21}+B_{22},s;\gamma)ds\\[0.2cm]
 &&+\gamma \int_{t_1}^{t_2}\bar{f}_{5}(I_0+\Xi_{11}+A_2,-\theta_{0} -\gamma^{-n}\psi_1-\Xi_{12}+\gamma^{-n}B_{21}+B_{22},s;\gamma)ds.\\[0.2cm]
\end{array}
$$
Again from (\ref{A1}), (\ref{A2}) and \eqref{xi1}, it holds that
$$|\Xi_{2i}(I_0,\theta_0; \gamma)|\leq O(\gamma^{5n+1})+\epsilon \,O(\gamma)$$
as $\gamma\to 0^+$ for $i=1,2$. As for the estimates on the derivatives, by equations (\ref{ABC2}),
we have
$$
\begin{array}{lll}
D_{I_0}B_{21}&=&\int_{t_1+}^tD_{I_0}\bar{f}_{3}\cdot\big(1+D_{I_0}(\Xi_{11}+A_2)\big)ds,\\[0.2cm]
D_{I_0}A_2&=&\gamma^{5n+1}\int_{t_1+}^t[D_{I_0}\bar{f}_{1}\cdot(1+D_{I_0}(\Xi_{11}+A_2))\\[0.2cm]
&&+D_{\theta_0}\bar{f}_{1}\cdot\big(\gamma^{-n}D_{I_0}(-\psi_1+B_{21})+D_{I_0}(-\Xi_{12}+B_{22})\big)]ds\\[0.2cm]
&&+\gamma \int_{t_1+}^t[D_{I_0}\bar{f}_{2}\cdot(1+D_{I_0}(\Xi_{11}+A_2))\\[0.2cm]
&&+D_{\theta_0}\bar{f}_{2}\cdot\big(\gamma^{-n}D_{I_0}(-\psi_1+B_{21})+D_{I_0}(-\Xi_{12}+B_{22})\big)]ds,\\[0.2cm]
D_{I_0}B_{22}&=&\gamma^{5n+1} \int_{t_1+}^t[D_{I_0}\bar{f}_{4}\cdot(1+D_{I_0}(\Xi_{11}+A_2))\\[0.2cm]
&&+D_{\theta_0}\bar{f}_{4}\cdot\big(\gamma^{-n}D_{I_0}(-\psi_1+B_{21})+D_{I_0}(-\Xi_{12}+B_{22})\big)]ds\\[0.2cm]
&&+\gamma \int_{t_1+}^t[D_{I_0}\bar{f}_{5}\cdot(1+D_{I_0}(\Xi_{11}+A_2))\\[0.2cm]
&&+D_{\theta_0}\bar{f}_{5}\cdot\big(\gamma^{-n}D_{I_0}(-\psi_1+B_{21})+D_{I_0}(-\Xi_{12}+B_{22})\big)]ds,
\end{array}
$$
which together with (\ref{A0}), (\ref{A1}), (\ref{A2}) and (\ref{xi1}) implies that
$$|D_{I_0}\Xi_{2i}(I_0,\theta_0; \gamma)|\leq O(\gamma^{4n+1})+\gamma^{-n}\epsilon \,O(\gamma)$$
for $i=1,2$. Differentiating $A_2,B_{21},B_{22}$ with respect to $I_0$ $k$ times yields that
$$|D_{I_0}^k\Xi_{2i}(I_0,\theta_0; \gamma)|\leq O(\gamma^{5n+1-kn})+\gamma^{-kn}\epsilon \,O(\gamma)$$
for $i=1,2$ and $0\leq k\leq 4$. Also by Differentiating $A_2,B_{21},B_{22}$ with respect to $\theta_0$, we know that
\begin{equation*}
\begin{array}{lll}
D_{\theta_0}B_{21}&=&\int_{t_1+}^tD_{I_0}\bar{f}_{3}\cdot D_{\theta_0}(\Xi_{11}+A_2)ds\\[0.2cm]
D_{\theta_0}A_2&=&\gamma^{5n+1} \int_{t_1+}^t[D_{I_0}\bar{f}_{1}\cdot D_{\theta_0}(\Xi_{11}+A_2)\\[0.2cm]
&&+D_{\theta_0}\bar{f}_{1}\cdot(-1-\gamma^{-n}D_{\theta_0}(\psi_1-B_{21})+D_{\theta_0}B_{22})]ds\\[0.2cm]
&&+\gamma \int_{t_1+}^t[D_{I_0}\bar{f}_{2}\cdot D_{\theta_0}(\Xi_{11}+A_2)\\[0.2cm]
&&+D_{\theta_0}\bar{f}_{2}\cdot(-1-\gamma^{-n}D_{\theta_0}(\psi_1-B_{21})+D_{\theta_0}B_{22})]ds,\\[0.2cm]
D_{\theta_0}B_{22}&=&\gamma^{5n+1} \int_{t_1+}^t[D_{I_0}\bar{f}_{4}\cdot D_{\theta_0}(\Xi_{11}+A_2)\\[0.2cm]
&&+D_{\theta_0}\bar{f}_{4}\cdot(-1-\gamma^{-n}D_{\theta_0}(\psi_1-B_{21})+D_{\theta_0}B_{22})]ds\\[0.2cm]
&&+\gamma \int_{t_1+}^t[D_{I_0}\bar{f}_{5}\cdot D_{\theta_0}(\Xi_{11}+A_2)\\[0.2cm]
&&+D_{\theta_0}\bar{f}_{5}\cdot(-1-\gamma^{-n}D_{\theta_0}(\psi_1-B_{21})+D_{\theta_0}B_{22})]ds,
\end{array}
\end{equation*}
which combining with (\ref{A0}), (\ref{A1}),(\ref{A2}) and (\ref{xi1}) implies that
$$|D_{\theta_0}\Xi_{2i}(I_0,\theta_0; \gamma)|\leq O(\gamma^{5n+1})+\epsilon \,O(\gamma), \ \ i=1,2.$$
The estimates on the higher derivatives $D_{I_0}^k D_{\theta_0}^l\Xi_{2i}$ can be proven similarly, from which we have
\begin{equation}\label{xi2}
|D_{I_{0}}^k D_{\theta_{0}}^l\Xi_{2i}(I_0,\theta_0; \gamma)|\leq O(\gamma^{n+1})+\gamma^{-4n}\epsilon \,O(\gamma)
\end{equation}
as $\gamma\to 0^+$ for $i=1,2$ and $k+l\leq 4$.

By impulsive functions in \eqref{I theta} and the value at $t_2-$ in \eqref{Pt2left}, we have
\begin{equation}\label{Pt2right}
\begin{array}{ll}
\quad \begin{array}{ll}
I(t_2+)=I(t_2-)=I_{0}+\Xi_{11}+\Xi_{21},\ \ \ \ \\[0.2cm]
\theta(t_2+)=-\theta(t_2-)=\theta_{0} +\gamma^{-n}\psi_1(I_{0};\gamma)-\gamma^{-n}\psi_2(I_{0};\gamma)+ \Xi_{12}- \Xi_{22},\ \ \\[0.2cm]
\end{array}
 \end{array}
 \end{equation}
With the same method, we finally consider $(I(t),\theta(t))$ on $[t_2,1]$ and obtain
$$
\begin{array}{lll}
I(t)&=&I(t_2+)+A_3(I_0,\theta_0,t;\gamma),\\[0.2cm]
\theta(t)&=&\theta(t_2+)+\gamma^{-n} B_{31}(I_0,\theta_0,t;\gamma)+B_{32}(I_0,\theta_0,t;\gamma),
\end{array}
$$
where functions $A_3,B_{31},B_{32}$ satisfy
\begin{equation}\label{ABC3}
\begin{array}{lll}
B_{31}&=&\int_{t_2}^t\bar{f}_{3}(I({t_2}+)+A_3,s;\gamma)ds\\[0.2cm]
A_3&=&\gamma^{5n+1} \int_{t_2}^t\bar{f}_{1}(I({t_2}+)+A_3,\theta({t_2}+)+\gamma^{-n}B_{31}+B_{32},s;\gamma)ds\\[0.2cm]
&&+\gamma \int_{t_2}^t\bar{f}_{2}(I({t_2}+)+A_3,\theta({t_2}+)+\gamma^{-n}B_{31}+B_{32},s;\gamma)ds,\\[0.2cm]
B_{32}&=&\gamma^{5n+1} \int_{t_2}^t\bar{f}_{4}(I({t_2}+)+A_3,\theta({t_2}+)+\gamma^{-n}B_{31}+B_{32},s;\gamma)ds\\[0.2cm]
&&+\gamma \int_{t_2}^t\bar{f}_{5}(I({t_2}+)+A_3,\theta({t_2}+)+\gamma^{-n}B_{31}+B_{32},s;\gamma)ds.
\end{array}
\end{equation}
Correspondingly, there are
\begin{equation}\label{P1}
\begin{array}{ll}
\quad \begin{array}{ll}
I(1)=I(t_2+)+\Xi_{31}(I_0,\theta_0; \gamma), \ \ \ \\[0.2cm]
\theta(1)=\theta(t_2+)+\gamma^{-n}\psi_3(I_{0};\gamma)+\Xi_{32}(I_0,\theta_0; \gamma),\ \ \\[0.2cm]
\end{array}
 \end{array}
 \end{equation}
where
\begin{equation}\label{psi3}
\begin{array}{lll}
\psi_3(I_0;\gamma)&=&\int_{t_2}^{1}\,\bar{f}_{3}(I(t),t;\gamma)dt\\[0.2cm]
&=&\int_{t_2}^{1}\,[d \cdot I(t)^n+ I(t)^{n-1}O(\gamma)]dt\\[0.2cm]
&=&(1-t_2)d \cdot \big(I(t_2+)+O(\gamma^{5n+1})+\epsilon O(\gamma)\big)^n\\[0.2cm]
&&+(1-t_2)\big(I(t_2+)+O(\gamma^{5n+1})+\epsilon O(\gamma)\big)^{n-1}O(\gamma )\\[0.2cm]
&=&(1-t_2)d \cdot I(t_2+)^n+(1-t_2)O(\gamma)(1+\epsilon)I(t_2+)^{n-1}+\cdots\\[0.2cm]
&=&(1-t_2)d \cdot \big(I_0+\Xi_{11}+\Xi_{21}\big)^n+(1-t_2)O(\gamma)(1+\epsilon)\big(I_0+\Xi_{11}+\Xi_{21}\big)^{n-1}+\cdots\\[0.2cm]
&=&(1-t_2)d \cdot \big(I_0+O(\gamma^{5n+1})+\epsilon O(\gamma)\big)^n\\[0.2cm]
&&+(1-t_2)O(\gamma)(1+\epsilon)\big(I_0+O(\gamma^{5n+1})+\epsilon O(\gamma)\big)^{n-1}+\cdots\\[0.2cm]
&=&(1-t_2)\big(d \cdot I_0^n+O(\gamma)(1+\epsilon)I_0^{n-1}+\cdots \big)
\end{array}
\end{equation}
and
$$
\begin{array}{lll}
 \Xi_{31}(I_0,\theta_0; \gamma)=A_3(I_0,\theta_0,1;\gamma)\\[0.2cm]
 \Xi_{32}(I_0,\theta_0; \gamma)=B_{32}(I_0,\theta_0,1;\gamma).
\end{array}
$$
Meanwhile, by the first two part discussion on $[0,t_1)$ and $[t_1,t_2)$, the estimates on the higher derivatives $D_{I_0}^k D_{\theta_0}^l\Xi_{3i}$ can be proven similarly, from which we have
\begin{equation}\label{xi3}
|D_{I_{0}}^k D_{\theta_{0}}^l\Xi_{3i}(I_0,\theta_0; \gamma)|\leq O(\gamma^{n+1})+\gamma^{-4n}\epsilon \,O(\gamma)
\end{equation}
as $\gamma\to 0^+$ for $i=1,2$ and $k+l\leq 4$.

Therefore, from \eqref{Pt2right} and \eqref{P1},
$$
\begin{array}{ll}
\quad \begin{array}{ll}
I(1)=I_{0}+\Xi_{11}+\Xi_{21}+\Xi_{31}, \ \ \ \\[0.2cm]
\theta(1)=\theta_{0} +\gamma^{-n}\big(\psi_1(I_{0};\gamma)-\psi_2(I_{0};\gamma)+\psi_3(I_{0};\gamma)\big)+ \Xi_{12}- \Xi_{22}+\Xi_{32}.\ \ \\[0.2cm]
\end{array}
 \end{array}
$$
Notice that by \eqref{psi1}, \eqref{psi2} and \eqref{psi3}
$$
\begin{array}{lll}
\psi_1(I_{0};\gamma)-\psi_2(I_{0};\gamma)+\psi_3(I_{0};\gamma)\\[0.2cm]
=[1-2(t_2-t_1)]\big(d \cdot I_0^n+O(\gamma)(1+\epsilon)I_0^{n-1}+\cdots \big)\\[0.2cm]
>\frac{d}{2}[1-2(t_2-t_1)],\ (\mbox{or}<\frac{d}{2}[1-2(t_2-t_1)])
\end{array}
$$
if $1-2(t_2-t_1)>0$ (or $1-2(t_2-t_1)<0$) for sufficiently small $\gamma$ and $I_{0}\in[1,2]$.
By denoting
$$
\Pi(I_{0};\gamma)\triangleq\psi_1(I_{0};\gamma)-\psi_2(I_{0};\gamma)+\psi_3(I_{0};\gamma),
$$
$$
\Xi_1( I_{0},\theta_{0}; \gamma)\triangleq\Xi_{11}+\Xi_{21}+\Xi_{31},
$$
and
$$
\Xi_2(I_{0},\theta_{0}; \gamma)\triangleq\Xi_{12}- \Xi_{22}+\Xi_{32},
$$
we rewrite $I(1)$ and $\theta(1)$ in the form of \eqref{estimation of P} and $\Pi$, $\Xi_1$, $\Xi_2$ satisfy Lemma \ref{estimation}. The proof is completed.\qed

\section{Proof of Theorems \ref{boundedness}- \ref{hamonic period solutions}}
In this section, we will conclude that the time 1 map $P$ is close to a twist map even though the existence of impulses. By means of a so-called large twist theorem established by X. Li, B. Liu and Y. Sun in \cite{XLi}, there exist large invariant curves diffeomorphic to circles and surrounding the origin in the $(x,x')$ plane. Every each curve is the base of a time periodic and under the flow invariant cylinder in the phase space $(x,x',t)\in\mathbb{R}^{2}\times\mathbb{R}$, which confines the solutions in its interior and which therefore leads to the boundedness of solutions. Finally, we construct the quasiperiodic solutions starting at $t=0$ on the invariant curves. Of course, since the existence of imlulses, the quasiperiodic solutions are discontinuous, which leads to the little different definition of corresponding shell functions. This is something new in this paper.

First we state the large twist theorem in \cite{XLi} which is used in the proof to obtain the invariant curves.
Consider the mapping
\begin{equation}\label{M}
\mathfrak{M}:\quad \begin{array}{ll}
\left\{\begin{array}{ll}
x_1=x+\beta +\gamma^{-\nu}\,y+ f(x,y;\gamma),\\[0.2cm]
y_1=y+g(x,y;\gamma),\\[0.1cm]
 \end{array}\right.\  \ \ \ \ (x,y)\in \mathbb{R} \times [a,b],
\end{array}
\end{equation}
where $b-a\geq 1$, $\nu\geq 0$ and $\beta$ are two constants, $0<\gamma\leq 1$ is a sufficiently small parameter. Here we assume that $f,g$ are real analytic in $x,y$, continuous in $\gamma$, and have period $2\pi$ in $x$, which can be extended to a complex domain
\begin{equation}\label{D1}
D:\ \  |\Im x|<r_0, \quad y\in D'
\end{equation}
with $D'$ a complex neighborhood of the interval $a\leq y\leq b$, $0<r_0<1$, and the mapping $\mathfrak{M}$ has the intersection property that any curve $\Gamma: y=\phi(x)=\phi(x+2\pi)$ always intersects its image curve $\mathfrak{M}\Gamma$.

We choose some $\omega$ satisfying
\begin{equation}\label{M1}
\beta+\gamma^{-\nu}\,a+\frac{1}{4}<\omega<\beta+\gamma^{-\nu}\,b-\frac{1}{4}
\end{equation}
and
\begin{equation}\label{M2}
\left|\frac{\omega}{2\pi}-\frac{q}{p}\right|\geq \frac{\gamma^\kappa}{p^\mu}
\end{equation}
for all integers $p,q$ with $p>0$ and $\mu>2$, $0<\kappa<\frac{1}{2}$. First we must show that for any sufficiently small $0<\gamma\leq 1$, there exists some $\omega$ satisfying (\ref{M1}) and (\ref{M2}). Then there is the large twist theorem.
\begin{theorem}\label{A}
Under these hypotheses, for each sufficiently small $\gamma>0$ and any number $\alpha\in(\frac{1}{2},1)$, there exists a positive constant $c$, depending on $D$, $\alpha$ and $\mu$ in (\ref{M2}), but not on $\gamma$, $\nu$, $\kappa$ and $\beta$, such that for
\begin{equation}\label{small}
\gamma^{\nu+\kappa}|f|+|g|<c\,\gamma^{\frac{\nu+\kappa}{1-\alpha}}
\end{equation}
in $D$ the mapping $\mathfrak{M}$ in (\ref{M}) admits an invariant curve of the form
\begin{equation}\label{Curve}
x=\xi+ p(\xi;\gamma), \ \ \ \ \ \ \ \ y=q(\xi;\gamma)
\end{equation}
with $p,q$ real analytic functions of period $2\pi$ in the complex domain $|\Im\, \xi|<\frac{r_0}{2}$ and continuous in $\gamma$. Moreover, the parametrization is chosen so that the induced mapping on the curve (\ref{Curve}) is given by
\begin{equation}\label{RN}
\xi_1=\xi+\omega
\end{equation}
with $\omega$ satisfying (\ref{M1}) and (\ref{M2}), and the functions $p,q$ satisfy
$$
\gamma^{\nu+1}|p|+|q-\gamma^{\nu}\,(\omega-\beta)|<\bar{c}\,\gamma^{\nu+\frac{\nu+\kappa}{1-\alpha}},
$$
where $\bar{c}$ is same as $c$, depending on $D$, $\alpha$ and $\mu$ in (\ref{M2}), but not on $\gamma$, $\nu$, $\kappa$ and $\beta$.
\end{theorem}

In order to guarantee the use of Theorem \ref{A}, we need the following essential intersection property of Poincar\'{e} map. In the left part, let
$$
\mathcal{A}\triangleq\{(I,\theta):\ I\in[1,2],\ \theta\in \mathbb{S}^1\}.
$$

\begin{lemma}\label{intersection}
The time 1 map $P$ has the intersection property on $\mathcal{A}$, i.e. if $C$ is an embedded circle in $\mathcal{A}$ homotopic to a circle $I=const$ in $\mathcal{A}$, then $P(C)\cap C\neq\emptyset$.
\end{lemma}
\Proof According to Section 4, $P$ is an area-preserving homeomorphism. Therefore, if $i_{C}:C\rightarrow \mathcal{A}$ is the injection map, we have
\begin{equation}\label{area preserving}
\displaystyle\int_{P(C)}i_{P(C)}^{\ast}\omega=\displaystyle\int_{C}i_{C}^{\ast}\omega,
\end{equation}
where $\omega=\lambda d\theta$.
Assume now, by contradiction, that $P(C)\cap C=\emptyset$, then $P(C)\cap C$ bounds an annulus in $\mathcal{A}$, which has positive measure. Consequently, by Stokes formula,
$$
\displaystyle\int_{P(C)}i_{P(C)}^{\ast}\omega-\displaystyle\int_{C}i_{C}^{\ast}\omega\neq 0,
$$
in contradiction to \eqref{area preserving}. This proves the Lemma.\qed

From Lemma \ref{estimation}, when $t_{2}-t_{1}=1/2$, $\Pi(I_{0};\gamma)$ in \eqref{estimation of P} equals to 0 and the time 1 map $P$ does not have a form of standard large twist map as \eqref{M} so that we cannot judge the boundedness of the solutions. Therefore we just consider the case $t_{2}-t_{1}\neq1/2$ in this paper. Since in the estimation \eqref{estimation of P} of $P$, there is a large perturbation term $\gamma^{-n}\Pi(I_{0};\gamma)$ in $\theta(1)$ as $\gamma\rightarrow 0$, which we need a large twist theorem in \cite{XLi} instead of Poincar\'{e} twist theorem.
Therefore, due to Lemma \ref{estimation}, we can choose $\epsilon=\gamma^{5n}$ with $\gamma$ small enough such that the functions $\Xi_1,\Xi_2$ satisfy the smallness condition in Theorem \ref{A}. The time 1 map $P$ is, with its derivatives, close to a large twist map. Moreover, it has intersection property in view of Lemma \ref{intersection}, such that the assumptions of the large twist theorem are met.

It follows that for $\omega\geq\omega^{\ast}$, $\omega^{\ast}$ sufficiently large, and
\begin{equation}\label{omega}
|\omega-\frac{p}{q}|\geq c|q|^{-2-\beta}
\end{equation}
for two constants $\beta>0$ and $c>0$ and for all integers $p$ and $q\neq 0$, there is an embedding $\psi:\mathbb{S}^{1}\rightarrow\mathcal{A}$ of a circle, which is invariant under the map $P$. Moreover, on this invariant curve the map $P$ is conjugated to a rotation with number $\omega$:
\begin{equation}\label{invariant map}
P\circ\psi(s)=\psi(s+\omega)\ \mbox{with}\ s(\ \mbox{mod}\ 1).
\end{equation}
Under the scalar transformation $\rho=\lambda^{\frac{1}{n+2}}$ and $\rho=\gamma^{-1}I$, there also exist invariant curves under $P$ satisfying \eqref{invariant map} in $(\lambda,\theta)$ coordinates. May as well denote the invariant curve by $\psi$ with rotation number $\omega$. Then the solutions of the Hamiltonian equation starting at time $t=0$ on this invariant curve determine a 1-periodic cylinder in the space $(\lambda,\theta, t)\in\mathcal{A}_{\lambda_{0}}\times\mathbb{R}$.
Since the Hamiltonian vector field $X_{H}$ \eqref{XH} is timeperiodic, the phase space is $\mathcal{A}_{\lambda_{0}}\times \mathbb{S}^{1}$. Let $\Phi^{t}$ with $\Phi^{0}=id$ be the flow of the time-independent vector field $(X_{H},1)$ on $\mathcal{A}_{\lambda_{0}}\times \mathbb{S}^{1}$ and define $\Psi(s,\tau):\mathbb{S}^{1}\times\mathbb{S}^{1}\rightarrow\mathcal{A}_{\lambda_{0}}\times \mathbb{S}^{1}$ by setting
$$
\Psi(s,\tau)=\Phi^{\tau}(\psi(s-\tau\omega),0)=(\phi^{\tau}\circ\psi(s-\tau\omega),\tau).
$$
Since $\phi^{t}$ is the flow of the vectorfield \eqref{I theta}, $\Psi(s,\tau)$ is discontinuous in $\tau$ actually.
In view of \eqref{invariant map}, we have $P=\phi^{1}$. Then
\begin{equation}\label{psi time period }
\begin{array}{ll}
\Psi(s,\tau+1)&=(\phi^{\tau+1}\circ\psi(s-(\tau+1)\omega),\tau+1)\\[0.2cm]
&=(\phi^{\tau}\circ\phi^{1}\circ\psi(s-\tau\omega-\omega),\tau+1)\\[0.2cm]
&=(\phi^{\tau}\circ\psi(s-\tau\omega),\tau)=\Psi(s,\tau),
\end{array}
\end{equation}
and
\begin{equation}\label{psi angle period }
\Psi(s+1,\tau)=\Phi^{\tau}(\psi(s+1-\tau\omega),0)=\Psi(s,\tau).
\end{equation}
Moreover,
\begin{equation}\label{psi quasi period }
\Phi^{t}\circ\Psi(s,\tau)=\Psi(s+\omega t,\tau+t)
\end{equation}
since $(X_{H},1)$ is an autonomous system.
Let
\begin{equation}\label{define quasi}
X(t;s,\tau)\triangleq\Psi(s+\omega t,\tau+t)
\end{equation}
with $X(0;s,\tau)=\Psi(s,\tau)$. By \eqref{psi quasi period }, there is
$$
X(t;s,\tau)=\Phi^{t}\circ X(0;s,\tau),
$$
which implies $X(t;s,\tau)$ is a solution of $(X_{H},1)$ starting from $X(0;s,\tau)$. Since $X(0;s,\tau)$ is the vector where the flow of $(X_{H},1)$ starting at $t=0$ on the invariant curve moves till $t=\tau$, combining with \eqref{psi time period } and \eqref{psi angle period } and choosing $(s,\tau)=(0,0)$ in \eqref{define quasi}, we conclude that solution $X(t;s,\tau)$ of $X_{H}$ is quasiperiodic with frequencies $(1,\omega)$ and $\Psi$ is its shell function, i.e.
$$
X(t;0,0)=\Psi(\omega t,t).
$$
Here we give some explanations about the quasiperiodic solutions. The discontinuous of $\Psi(s,\tau)$ results that $X(t)$ is piecewise continuous. On each interval without impulses, $X(t)$ has the form of quasiperiodic solutions while at impulsive times $X(t)$ satisfies impulsive conditions. Strictly speaking, $\Psi(\omega t,t)$ is not a continuous shell function for all $t\in\mathbb{R}$. But according to the form of $X(t)$ at $t\neq t_{j},\ j=\pm1,\pm2,\cdots$, we also regard $X(t)$ is quasiperiodic with piecewise continuous shell function $\Psi$.

Through the transformation \eqref{symplectic tran}, solutions of $X_{h}$ \eqref{Xhfull} are quasiperiodic either. This prove the statement of Theorem \ref{quasiperiodic solution}. In order to prove the statement of Theorem \ref{boundedness}, we transform the invariant curves obtained for $X_{H}$ into the $(x,y)$ coordinates. Since the transformation $\varphi$ and $\psi$ are symplectic, there are invariant curves of the time 1 map of the flow which go around the origin and extend to infinity. Its solution is therefore confined in the interior of the time periodic cylinder above the invariant curve and hence is bounded. This ends the proof of Theorem \ref{boundedness}. \qed

We next prove Theorem \ref{hamonic period solutions}. Following the argument of G. Morris \cite{morris} we shall establish fixed points of the iterated map $P^{m}$ for every integer $m$ applying the Poincar\'{e}-Birkhoff fixed point theorem to the map $P^{m}$ in the annulus $R$ bounded by two invariant curves with rotation numbers $\omega_{1}=\omega<\omega_{2}=\omega+1$ satisfying \eqref{omega}. First, we map this annulus $R$ onto an annulus $R_{0}=\{(\xi,\eta)|\ 0\leq\eta\leq 1,\xi\ \mbox{mod}\ 1\}$ bounded by concentric circles. If $\psi_{j}$ are the embeddings of the invariant curves (see \eqref{invariant map}) having rotation numbers $\omega_{j}$, $j=1,\ 2$, we define the map $\chi:R_{0}\rightarrow R$ by
$$
\chi:(\mu, \theta)=\eta\psi_{1}(\xi)+(1-\eta)\psi_{2}(\xi).
$$
Since the embeddings $\psi_{j}$ are differentiably close to the injections of the circles $\omega_{j}=const$, $\chi$ is a diffeomorphism. The induced map $G:\chi^{-1}\circ P\circ\chi$, expressed by
\begin{equation}\label{G}
\begin{array}{ll}
\widehat{G}:\quad \begin{array}{ll}
\xi_{1}=f(\xi,\eta),\\[0.2cm]
\eta_{1}=g(\xi,\eta).
\end{array}
\end{array}
\end{equation}
satisfies in view of \eqref{invariant map} $f(\xi+1,\eta)=f(\xi,\eta)$ and $g(\xi+1,\eta)=g(\xi,\eta)$. It leaves the boundaries $\eta=0$ and $\eta=1$ invariant and satisfies on these boundaries:
\begin{equation}\label{0 1 bondary}
f(\xi,0)=\xi+\omega_{2},\ \ \ \ f(\xi,1)=\xi+\omega_{1}.
\end{equation}
Define the map $s:R_{0}\rightarrow R_{0}$ by $s(\xi,\eta)=(\xi+1,\eta)$. To find a period point of minimal period $q\geq1$ of $G$, for any given integer $q$, we let $[q\omega]$ be the integral part of $q\omega$, i.e. $q\omega-1\leq[q\omega]\leq q\omega$, and choose an integer $0<p\leq q$ such that $p+[q\omega]$ and $q$ are relatively prime. The map $G_{q}:=s^{-p-[q\omega]}\circ G^{q}$ on $R_{0}$, expressed by
\begin{equation}\label{Gq}
\begin{array}{ll}
G_{q}:\quad \begin{array}{ll}
\xi_{1}=f_{q}(\xi,\eta),\\[0.2cm]
\eta_{1}=g_{q}(\xi,\eta).
\end{array}
\end{array}
\end{equation}
satisfies the required twist conditions at the boundaries. Indeed with $0<\alpha=q\omega-[q\omega]<1$ we conclude from \eqref{0 1 bondary}:
$$
f_{q}(\xi,1)-\xi=-p+\alpha<0<(q-p)+\alpha=f_{q}(\xi,0)-\xi.
$$
Since $G_{q}$ leaves a regular measure invariant we conclude that $G_{q}$ possesses a fixed point in $R_{0}$. One verifies easily that it corresponds to a periodic point of $G$, hence of $P$, with minimal period $q$. This proves the statement of Theorem \ref{hamonic period solutions}.\qed

\section*{References}
\bibliographystyle{elsarticle-num}

\end{document}